\documentclass[review,3p,a4paper]{elsarticle}

\usepackage{amssymb}
\usepackage{amsthm}
\usepackage{amsmath,bbm,dsfont, color}

\journal{arXiv}
\begin{document}
\renewcommand{\theequation}{\thesection.\arabic{equation}}
\overfullrule5pt
\newcommand*{\id}{\mathds{1}}
\newcommand*{\zahlen}{\mathbb{Z}}
\newcommand*{\en}{\mathbb{N}}
\newcommand*{\er}{\mathbb{R}}
\newcommand{\esssup}{\mathop{\text{ess\:sup}}}
\newcommand{\htimes}{\mathop{\text{\large$�$}}}
\newcommand{\hexists}{\mathop{\text{\LARGE$\exists$}}}
\newcommand{\hforall}{\mathop{\text{\LARGE$\forall$}}}
\newtheorem{theo}{Theorem}[section]                    
\newtheorem{lem}[theo]{Lemma}
\newtheorem{cor}[theo]{Corollary}
\newtheorem{prop}{Proposition}[section]
\newtheorem{mydef}{Definition}[section]
\newtheorem{rem}{Remark}[section]
\newcommand{\Sym}{{\rm Sym}}

 \def\Xint#1{\mathchoice
 {\XXint\displaystyle\textstyle{#1}}%
 {\XXint\textstyle\scriptstyle{#1}}%
 {\XXint\scriptstyle\scriptscriptstyle{#1}}%
 {\XXint\scriptscriptstyle\scriptscriptstyle{#1}}%
 \!\int}
 \def\XXint#1#2#3{{\setbox0=\hbox{$#1{#2#3}{\int}$}
 \vcenter{\hbox{$#2#3$}}\kern-.5\wd0}}
 \def\ddashint{\Xint=}
 \def\dashint{\Xint-}

\def\divv{{\rm div}}
\def\B{\mathbb{B}}
\def\D{\mathbb{D}}
\def\F{\mathbb{F}}
\def\N{\mathbb{N}}
\def\O{\mathbb{O}}
\def\R{\mathbb{R}}
\def\Kor{{\rm Kor}}
\def\pod#1{\mathop{#1}\limits}
\def\diagin{-\hskip-11.0truept\intop}
\def\diagint{{\raise-.1pt\hbox{--}\hskip-7.9pt\intop}}
\def\diagintop{\mathop{\mathchoice
{{\diagin}}%
{{\diagint}}%
{{\diagint}}%
{{\diagint}}%
}\limits}

\begin{frontmatter}

 \title{Almost everywhere H\"older continuity of gradients to non-diagonal parabolic systems}

\author{Jan Burczak} \ead{jb@impan.pl} \address{Institute of Mathematics, Polish Academy of Sciences, \'Sniadeckich 8, 00-950 Warsaw.}






\begin{abstract}
We present a local almost everywhere regularity result for a general p-nonlinear non-diagonal parabolic system, the main part of which depends on symmetric part of the gradient.
\end{abstract}





\end{frontmatter}

\section{Introduction}

\noindent
The problem of local H\"older continuity of gradients for the evolutionary $p$-Laplace system
has been resolved in a series of papers by DiBenedetto and coauthors, summed up in a monograph \cite{[B]}, with crucial earlier (stationary) contributions of Uhlenbeck \cite{Uhl}, Tolksdorf \cite{[T]} and the Russian school. From the perspective of mathematical physics, it is interesting to replace $\nabla u$ by its symmetric part $\D u=(\nabla u+\nabla^Tu) / 2$; then such a symmetric $p$-Laplace system is a simplification of the hydrodynamic model of a non-Newtonian flow (referred to as $p$-Navier-Stokes in the following). In fact, for $ p > 11/5$, the generalization from a $p$-Stokes system to the respective hydrodynamic one is not essential from the perspective of regularity theory (compare \cite{[NM]}). \\
 It turns out that the amendment from $\nabla u$ to $\D u$ in the $p$-Laplace system, a supposedly harmless one, diminishes dramatically our understanding of $C^{1,\alpha}$-regularity of such system. The reason is that most of the relevant methods successful in the full gradient case turn out to be useless, because they rely essentially on pointwise structure. In this article we show, however, that the caloric approximation approach can still be used to obtain almost everywhere regularity.
We consider parabolic systems of the following type
\begin{equation}\label{1.1}
u_{,t}-\divv A(z,u,\D u)=0
\end{equation} 
the prototype of which is the following {\it symmetric $p$-Laplace} system with safety $1$
\begin{equation}\label{1.1prot}
u_{,t}-\divv \left[(1 + |\D u|^2)^{\frac{p-2}{2}} \D u \right]=0
\end{equation}
Let us provide the reader with a short account of relevant known results. In \cite{[PB]} an extensive short-time maximal regularity theory in Sobolev-Slobodeckii spaces for $p$-Navier-Stokes is presented by Pr\"uss and Bothe. 
However, not much is known on the global-in-time $C^{1,\alpha}$-regularity of such systems in arbitrary dimension $d$ (or at least for physically plausible $d \ge 3$), even for the prototype case. 
It is worth mentioning that for $p \in (12/5; 10/3)$ Seregin has shown in \cite{[SL]} an almost-everywhere regularity result for the complete three-dimensional hydrodynamic system. One can also easily see from the theory developed in \cite{[KMS]} by Kaplick\'y, M\'alek, Star\'a for the two-dimensional $p$-Navier-Stokes that system \eqref{1.1prot} and some of its generalizations enjoy  $C^{1,\alpha}$-regularity  in the case of two-dimensions. The current research status is a little clearer in the case of stationary simplifications. There is a well developed $C^{1,\alpha}$-theory for the stationary $p$-Laplace and $p$-Navier-Stokes systems with $p<2$ by Bair\~ao da Veiga and collaborators (see \cite{[BdV]} and references therein). In the case of $p \ge 2$,  one can refer to \cite{[ABF]}, where Apushkinskaya, Bildhauer and Fuchs obtain partial $C^{1,\alpha}$-regularity for  three-dimensional $p$-Stokes and full $C^{1,\alpha}$-regularity in the two-dimensional case.\\
In this paper we follow the theory based on $p$-caloric approximations, which has been developed for the full-gradient case and very general main parts in \cite{[DM]}, \cite{[DMS]} by Duzaar, Mingione and coauthors. We apply their ideas for the symmetric-gradient case. At some points we could have merely quoted the respective results from \cite{[DMS]}; instead,  for reader's convenience, most of the proofs are presented with concern for the clarity of exposition. \\
However the $p$-caloric approach seems to be very well-suited also for our symmetric-gradient case, let us emphasize that our result seems to be new not only for general system \eqref{1.1}, but even for its prototype \eqref{1.1prot}. As a byproduct, we obtain also a Campanato-type theory for linear parabolic systems satisfying Legendre-Hadamard conditions, for which we couldn't find a satisfactory reference.
\section{Notation and statement of the result}
\setcounter{equation}{0}
The expression $A \equiv B$ means that $A$ is defined as $B$. Denote a space-time point $z=(x,t)\in \Omega \times(-T,0) \equiv Q$, where $\Omega \subset \er^d$. As we develop a local interior regularity theory, any further assumptions on domain $Q$ are unnecessary. $B_r (x), \; Q_r (z)$ denote, respectively, the ball with the radius $r$ centered at a point $x$ and the parabolic cylinder $B_r (x) \times (t-r^2, t)$. $\partial_\Gamma Q$ denotes parabolic boundary of cylinder $Q$.\\
For a tensor $\xi\in\R^{d\times d}$ denote its symmetric part by $\xi^s \equiv {(\xi+\xi^T) / 2}$. For any matrix $M \in \er^{d^2 \times {d^2}}$ denote it coefficients by $M^{ij}_{kl}$; its action on tensor  $\xi$ with coefficients $\xi_{kl}$ is  $M^{ij}_{kl} \xi_{kl}$ (here and further on we use the summation convention). $Sym^{d \times d}$ denotes set of $d \times d$ symmetric tensors.\\
We use standard notation for function spaces; $L^p (\tau, t; W^{1,p}(B_\varrho  (z) ) )$ will be sometimes abbreviated to $L^p ( W^{1,p})$, when there is no danger of confusion regarding underlying cylinder.\\
Let us emphasize that constants denoted by $C$ may change from line to line of estimates and are generally bigger than $1$. If a more careful control over a constant is needed, we denote their dependence on certain parameters writing $C(parameter)$ and generally suppress marking their dependence on irrelevant parameters; such constants may also vary. For clarity we also use some fixed constants, which we denote by $C_{subscript}$.\\
Now let us present a list of assumptions for the studied generalization \eqref{1.1} to  \eqref{1.1prot}. For any  tensors $\xi, \eta \in \R^{d\times d}$
\begin{itemize} \item[$\bullet$] main part $A$ satisfies properties of
\begin{equation}\label{1.2a}
A(z,u,\xi^s)\xi^s\ge\lambda|\xi^s|^p \qquad \text{(being {\it strongly elliptic})},
\end{equation}
\begin{equation}\label{1.3}
A(z,u,\xi^s)\eta \ge A(z,u,\xi^s)\eta^s \qquad \text{(being {\it weakly symmetrizing})},
\end{equation}
\begin{equation}\label{1.6}
|A(z,u,q)|\le C \cdot \left(1+|q|^{p-1} \right)  \qquad \text{(having $p-1$ growth)},
\end{equation}
\begin{equation}\label{1.8}
|A(z,u,q)-A(\tilde z,\tilde u,  q)|\le C \cdot \min \left(1, K(|u|+|\tilde u|)  \cdot \left( d_2 ( z - \tilde z) +|u-\tilde u|)^\beta \right)  \cdot (1+|q|^{p-1} \right),
\end{equation}
where $\beta\in(0,1)$ and $K: [0, \infty) \to [1, \infty)$ is a non-decreasing real function;

\item[$\bullet$] whereas ${\partial A\over \partial q}$
 \begin{equation}\label{1.2b}
{\partial A\over \partial q}(z,u,q)\xi^s\cdot\xi^s\ge \lambda (1 + |q|^2)^{\frac{p-2}{2}} |\xi^s|^2 \qquad \text{(is { \it Legendre-Hadamard elliptic})},
\end{equation} 
\begin{equation}\label{1.4}
\left({\partial A\over\partial q}(z,u,q) \right)^{ij}_{kl} = \left({\partial A\over\partial q}(z,u,q) \right)^{kl}_{ij}  =\left({\partial A\over\partial q}(z,u,q) \right)^{ji}_{lk} \qquad \text{(is { \it strongly symmetrizing})},
\end{equation}
\begin{equation}\label{1.7}
|u|+|q| \le M \; \implies \; \left|{\partial A\over\partial q}(z,u,q)\right|\le C_\eqref{1.7}(M) \qquad \text{(grows in a general way)},
\end{equation}
\begin{equation}\label{1.5}
|u|+|q|+|u-\tilde u|+|q-\tilde q|\le M \implies \left|{\partial A\over\partial q}(z,u,q)-{\partial A\over\partial \tilde q}(\tilde z,\tilde u,\tilde q)\right| \le C (M) \;\omega(M, d_2^2 (z-\tilde z)+ |u-\tilde u|^p+|q-\tilde q|^p)  \quad \text{(is continuous)}
\end{equation}
with parabolic metric $d_2(z-z_0) = |x - x_0| + |t - t_0|^\frac{1}{2}$
and local modulus of continuity $\omega$ satisfying: $\omega(\cdot,s)$, $\omega(t,\cdot)$ are nondecreasing, $\omega(t,0)=0$ and $\omega(t, \cdot)$ is continuous at zero, $\omega^{p}(t,\cdot)$ is concave.
\end{itemize}
\begin{rem}
Observe that property \eqref{1.5} is indeed merely continuity and that  \eqref{1.4} implies that ${\partial A\over \partial q}$ is weakly symmetrizing, i.e.
\begin{equation}\label{1.4'}
{\partial A\over\partial q}(z,u,q)\xi^s\cdot\eta=
{\partial A\over\partial q}(z,u,q)\xi^s\eta^s
\end{equation}
\end{rem}
The main result reads.
\begin{theo}\label{theo1.1} 
Any weak solution $u\in C(-T, 0; L^2 (\Omega) ) \cap L^p(-T, 0;  W^{1,p} (\Omega) )$ to the system (\ref{1.1}) with $p \ge 2$ and structure (\ref{1.2a} ---\ref{1.5}) has a.e. H\"older continuous gradients and the solution itself is also a.e. H\"older continuous. More precisely, there is an open set $\tilde Q$ of full Lebesgue measure satisfying
\begin{equation*}
\tilde Q \quad \supset \quad \left\{z \in Q:  \; \liminf_{\varrho\to0} \diagintop_{Q_\varrho(z)} |\nabla u-\ (\nabla u)_{z}  |^p =0 \quad \wedge \quad \limsup_{\varrho\to0}|(u)_{z,\varrho}|+|(\nabla u)_{z,\varrho}|<+\infty \right\} 
\end{equation*}
for which

\begin{equation*}
\nabla u\in C^{\beta,{\beta\over 2}}(\tilde Q), \quad  u\in C^{1,{1\over 2}}(\tilde Q),
\end{equation*}
where $\beta$ comes from (\ref{1.8}).

\end{theo}

Let us repeat that, to our best knowledge, even for the prototype system \eqref{1.1prot} the results is new.

\section{Outline of the paper}
The rest of the article is devoted to the proof of the result stated above. For traceability, let us first present the outline of the paper. In Section \ref{sec2} auxiliary lemmas are gathered. This includes a Campanato-type regularity theory for linear parabolic systems satisfying Legendre-Hadamard conditions, see Lemma \ref{LinSys} and the {\it symmetric caloric approximation lemma} --- Lemma \ref{lem5.1}. The latter states, in the context of symmetric gradients, that every function which is close to a solution of a linear parabolic system in a certain weak sense is indeed close to a solution of a linear parabolic system in a strong sense. Next sections are devoted to the proof of Theorem \ref{theo1.1}, the main steps of which are as follows.
\begin{enumerate}
\item Section \ref{sec:ineq} is devoted to showing, by means of linearization and Caccioppoli inequality, that an appropriately rescaled weak solution to \eqref{1.1} satisfies locally certain inequalities that resemble assumptions of the caloric approximation lemma. This is done via Lemmas  \ref{lem3.1},  \ref{lem4.1} and summed up in Corollary \ref{cor4.2}.
\item Section \ref{partial} combines results of the previous sections and gives the proof of Theorem \ref{theo1.1}. Namely, thanks to Corollary \ref{cor4.2} around points which satisfy certain regularity assumptions one can use caloric approximation for (rescaled) solution of \eqref{1.1}, which thanks to the regularity of linear systems gives proper shrinking of excess energies (Lemma \ref{lem6.1}). This yields, by iteration, the H\"older continuity of gradients (Lemma \ref{lem6.2}). Finally, the full thesis of the main theorem is obtained by  redoing estimates of previous Lemmas at the level of solutions (Lemma \ref{lem6.2'}).
\end{enumerate}
Only the crucial results are proved directly after their statements; for the sake of clarity, the remaining proofs are transferred to the Section \ref{app} --- Appendix.
\section{Useful auxiliary results} \label{sec2}
\setcounter{equation}{0}
This section begins with a Simon-type compactness result for parabolic spaces, which can be found as Theorem 2.5 in \cite{[DMS]}.
\begin{lem}\label{lem1.15}
Take $p \in (1, \infty)$, three Banach spaces $X \subset \subset Y \subset Z$ and a sequence $g_k$, which is uniformly bounded in $L^p(-T, 0 ; X)$ and satisfies
\begin{equation}\label{sim}
\hforall_{\varepsilon > 0} \hexists_{h'} \hforall_{h \in (0, h']} \intop_{-T}^{-h}|  g_k (\cdot,t+h)-  g_k(\cdot,t)|^p_{Z} dt \le \varepsilon
\end{equation}
then $g_k$ contains a subsequence convergent in the space $L^p(-T, 0 ; Y)$.
\end{lem}
The next result collects properties needed to perform analysis of excess energies. For proof see \cite{[DMS]} Lemma 2.1; the last inequality can be found in proof of Lemma 4.8 there. Compare also \cite{[K]}.
\begin{lem}\label{lem1.2}
Let $u\in L^s(Q(z_0))$, $s \ge 2$, where $z_0 = (x_0, t_0)$. There is the unique minimizer $l^{(s)}_\varrho(x)$ to $\intop_{Q_\varrho(z_0)}|u-l|^s$ among affine, time-independent functions $l$; moreover
\begin{equation}\label{1.2.1}
l^{(2)}_\varrho(x)=(u)_{z_0,\varrho}+\underbrace{ \left[ {d+2\over\varrho^2}\diagintop_{Q_\varrho(z_0)}u(x,t)\otimes(x-x_0) dx dt\right]}
_{Q^{(2)}_{z_0,\varrho}}(x-x_0)
\end{equation}
the linear part $Q^{(2)}_{z_0,\varrho}$ of which is close to $(\nabla u)_{z_0,\varrho}$
\begin{equation}\label{1.2.2}
|Q^{(2)}_{z_0,\varrho}-(\nabla u)_{z_0,\varrho}|^2\le {d(d+2)\over \varrho^2}   \diagintop_{Q_\varrho(z_0)}| u-(u)_{z_0,\varrho} - (\nabla u)_{z_0,\varrho} (x-x_0) |^2
\end{equation}
and shrinks as follows
\begin{equation}\label{1.2.3}
|Q^{(2)}_{z_0,\theta\varrho}-Q^{(2)}_{z_0,\varrho}|^2\le{d(d+2)\over(\theta\varrho)^2}\diagintop_{Q_\varrho(z_0)}|u-l^{(2)}_\varrho|^2
\end{equation}
For the minimizer in the case of general $s \ge 2$ holds
\begin{equation}\label{1.2.4}
\diagintop_{Q_\varrho(z_0)} \left|l^{(2)}_{{\varrho}}-l^{(s)}_{{\varrho}}\right|^s \le C_\eqref{1.2.4} (d,s) \diagintop_{Q_\varrho(z_0)} \left|u-l^{(s)}_{{\varrho}}\right|^s 
\end{equation}
\end{lem}
Subsequently let us state the Korn's inequality. For hints for proof, see the Appendix.
\begin{lem}\label{lem1.3} 
(Korn's inequality) 
For $u\in W^{1,p}(B_r(x))$ following inequalities hold with $K_p$ independent on radius of $B_r(x)$
\begin{equation}\label{1.3.1}
K_p\left[\intop_{B_r(x)} r^{-p} |u|^p+|\D u|^p\right]\ge\intop_{B_r(x)}|\nabla u|^p
\end{equation}
\begin{equation}\label{1.3.2}
K\intop_{B_r(x)}|\D u-(\D u)|^2\ge\intop_{B_r(x)}|\nabla u-(\nabla u)|^2\ge\intop_{B_r(x)}
|\D u-(\D u)|^2
\end{equation}
\end{lem}
Next lemma, which may be of independent interest,  collects needed results on linear parabolic systems with main part depending on symmetric gradient. Recall that $AM$ denotes constant coefficient matrix $A$ with elements  $a_{kl}^{ij}$ acting on tensor $M$ with elements  $m_{kl}$, i.e. $AM=a_{kl}^{ij}m_{kl}$. Again we refer to the Appendix for the proof.
\begin{lem}\label{LinSys}
(Campanato-type regularity theory for linear parabolic systems satisfying Legendre-Hadamard conditions) 
Let $u\in L^2(-T,0; W^{1,2} (\Omega))$ be a local solution to $u_{,t}-\divv A \; \D u=0$, i.e let it satisfy 
\begin{equation}\label{linC}
\intop_{\Omega_T}u\varphi_{,t}- A \; \D u \; \D\phi =0 \quad \hforall_{\varphi\in C_0^\infty(\Omega_T)}
\end{equation}
where for constant coefficient matrix $A$ holds:
\begin{equation}\label{2.5.1}
a_{kl}^{ij}=a_{ij}^{kl}=a_{lk}^{ji}
\end{equation}
\begin{equation}\label{2.5.2}
a_{kl}^{ij}\xi_l\xi_j\eta^k\eta^i+a_{kl}^{ij}\xi_k\xi_j\eta^l\eta^i\ge \lambda|\eta|^2|\xi|^2\ 
\hforall_{\eta,\xi\in\R^d}
\end{equation}
\begin{equation}\label{2.5.3}
A\xi^s\xi^s\ge\lambda|\xi^s|^2\ \hforall_{\xi\in\R^{d\times d}}
\end{equation}
then $u$ is locally smooth and satisfies for any $p, q \in [1, \infty], \; \varrho \le r/2$, arbitrary $\tilde{z}_0 \in Q_{\varrho}$
\begin{equation}\label{2.4.1}
\left[ \dashint_{Q_{\varrho}} \left| u^{(m)} \right|^q   \right]^\frac{1}{q} \le C_\eqref{2.4.1} (\lambda, |A|, K_p, m, d, p, q) \; r^{-2m}  \left[  \dashint_{Q_{r}} |u|^p \right]^{1/p},
\end{equation}
\begin{equation}\label{2.4.11}
\left[ \dashint_{Q_{\varrho}} \left|u^{(m)}-u^{(m)}(\tilde{z}_0) \right|^q \right]^\frac{1}{q}\le  C_\eqref{2.4.11} (\lambda, |A|, K_p, m, d, p, q) \; r^{-2m} \left(\frac{\varrho}{r} \right)  \left[  \dashint_{Q_{r}} |u|^p \right]^{1/p},
\end{equation}
\begin{equation}\label{2.4.12}
\left[  \dashint_{Q_{\varrho}} \left|u^{(m)}-\left(u^{(m)}\right)_{Q_{\varrho}} \right|^q \right]^\frac{1}{q} \le  C_\eqref{2.4.12} (\lambda, |A|, K_p, m, d, p, q) \; r^{-2m}  \left(\frac{\varrho}{r} \right)  \left[  \dashint_{Q_{r}} |u|^p \right]^{1/p}.
\end{equation}
where $u^{(m)}$ denotes either $\nabla^{(2m)} u$ or $\partial_t^{(m)} u$ and $ | a |^s  = \sum_{n = 1}^N | a_i |^s $ for $a \in \er^N$ .
\end{lem}
As outlined in the introduction, we end this section by stating a local result which says that a function, which is approximately solving a certain linear system in a weak sense (such function is called {\it $\delta$-approximatively weakly symmetrical caloric} in the following), is indeed close to some solution to this system in an appropriate strong $L^2 - L^p$ sense. The idea can be traced back to L. Simon, see \cite{[Sim]}. The proof, up to few technicalities connected with symmetric gradient, is identical with its counterpart in \cite{[DMS]} and can be found in the Appendix.
We work now with fixed $p \ge 2$ and cylinder $Q_\varrho (z_0)$ (therefore they does not appear as parameters). Let us introduce some definitions.
\begin{mydef}\label{def5.05}
$S(\lambda, \Lambda)$ denotes the set of elliptic bilinear forms, which have the properties of being symmetrizing and $\lambda$-elliptic and $\Lambda$-bounded. Precisely:
\begin{equation}
S(\lambda, \Lambda) := \left\{A\colon\R^{d^2}\times\R^{d^2}\to\R, \textsl{bilinear}, \quad
a_{kl}^{ij}=a_{lk}^{ji},\quad \lambda|\xi^s|^2 \le A\xi^s\xi^s,\quad  |A| \le \Lambda \quad
\hforall_{\xi,\eta\in\R^{d^2}} \right\}
\end{equation}
\end{mydef}
Observe that A is sweakly ymmetrizing, as $a_{kl}^{ij}=a_{lk}^{ji}$ implies $A\xi^s\eta=A\xi^s\eta^s$.\\
In the following two definitions $\delta >0, \gamma \ge 0$ are number parameters.
\begin{mydef}\label{def5.1}
Set $H (r; \delta, A, \gamma)$ of  approximatively weakly symmetrical caloric functions consists of elements of \\ $L^p (t_0 - r^2, t_0; W^{1,p}(B_r  (z_0) )$~that~satisfy
\begin{equation}\label{5.1.1}
\diagintop_{Q_r  (z_0)} \left|\frac{f}{r} \right|^2+ |\nabla f|^2 + \gamma^{p-2} \left[ \left|\frac{f}{r} \right|^p+ |\nabla f|^p \right] \le 1,\; \; 
\Bigg|\diagintop_{Q_r  (z_0)}f\varphi_{,t}-A(\D f,\D\varphi)\Bigg|\le\delta\cdot|\D\varphi|_{L^\infty} \quad
\hforall_{\varphi\in C_0^\infty({Q_r  (z_0)})}
\end{equation}
\end{mydef}
\begin{mydef}\label{def5.2}
Set $H (r; A, \gamma)$ of caloric symmetrical functions constitute  $f \in L^p (t_0 - r^2, t_0; W^{1,p}(B_r  (z_0) ) $ such that
\begin{equation}
\diagintop_{Q_r  (z_0)} \left|\frac{f}{r} \right|^2+ |\nabla f|^2 + \gamma^{p-2} \left[ \left|\frac{f}{r} \right|^p+ |\nabla f|^p \right] \le 2^{d+3}, \; \;
\diagintop_{Q_r  (z_0)}f\varphi_{,t} - A(\D f,\D\varphi) = 0 \quad
\hforall_{\varphi\in C_0^\infty({Q_r  (z_0)})}
\end{equation}
\end{mydef}
\begin{lem}[symmetric caloric approximation lemma]\label{lem5.1}
Take $p \ge 2$. Fix positive $\varepsilon, \lambda, \Lambda$. Then there exists $\delta \in (0,1)$, common for: all ${A \in S(\lambda, \Lambda)}$ and ${\gamma \in [0,1]} $, such that the following implication holds
\begin{equation*}
f\in H(\varrho; \delta, A, \gamma) \Rightarrow \hexists_{h\in H(\varrho/2;  A, \gamma)} \; \diagintop_{Q_{\varrho/2} (z_0)}
\left|\frac{h-f}{\varrho/2}\right|^2 + \gamma^{p-2} \left|\frac{h-f}{\varrho/2}\right|^p  \le \varepsilon
\end{equation*}
\end{lem}
\section{Local estimates}\label{sec:ineq}
\setcounter{equation}{0}
Let us emphasize that in this section the dependence of constants $C$ on irrelevant parameters is suppressed. First let us define local excess energies
\begin{mydef}\label{inxs}
\begin{equation}
\phi_{p, z_0, l} (\varrho)=\diagintop_{Q_{z_0}(\varrho)}|\D u-\D l|^p, \quad \psi_{p, z_0, l}   (\varrho) =\diagintop_{Q_{z_0}(\varrho)} \left|{u-l\over\varrho}\right|^p
\end{equation}
\end{mydef}
For briefness, using the energies defined above we often drop certain parameters, writing for example $\phi_p (\varrho),  \psi_{p} (\varrho)$. First we state an auxiliary algebraic lemma needed for the estimates of this section.
\begin{lem}[Algebraic inequalities]\label{lem3.05} Fix $M$. Assume that for matrix $A$ condition \eqref{1.8} is valid. Then, for any $z \in Q_{z_0} (\varrho) \subset Q$ with $\varrho \le 1$,
any $u \in \er^d, P \in Sym^{d \times d}$ and any affine function $l(x) $ the following inequalities hold 
\begin{equation}\label{alg_1}
|A(z_0, l(z_0), P) - A(z, u, P)| \le C(|l(z_0)| + |\nabla l| ) \varrho^\beta
\left [1+  |P-\D l|^{p} + \left| \frac{u - l}{\varrho} \right|^p \right] 
\end{equation}
\begin{equation}\label{alg_3}
|A(z_0,l(z_0), \D l) - A(z, u, \D l)| \le C(|l(z_0)| + |\nabla l| ) \varrho^\beta \left[1+ \left| \frac{u - l}{\varrho} \right|^\beta \right] 
\end{equation}

If, additionally, $A$ satisfies \eqref{1.6}, \eqref{1.7}, then it also holds
\begin{equation}\label{alg_2}
|A(z, u, P) - A(z, u, \D l)| \le C (|l(z_0)| + |\nabla l|)  \left( |u-l|^\beta+ |u-l| + |P-\D l|^{p-1}  \right) 
\end{equation}
\begin{equation}\label{alg_4}
|A(z, u, P) - A(z_0,l(z_0), \D l)| \le C(|l(z_0)| + |\nabla l|) (1 + |P|^{p-1}  )
\end{equation}
\end{lem}
Proof of this lemma has been shifted to appendix.

\begin{lem}[Linearization]\label{lem3.1} 
Take $p\ge2$ and $u \in L^p(-T, 0; W^{1,p} (\Omega))$ solving (\ref{1.1}) with structure: (\ref{1.3}),  (\ref{1.6}), (\ref{1.8}), (\ref{1.7}), (\ref{1.5}). For any $M>0$, $Q_\varrho(z_0)\subset Q$ with $\varrho \le 1$, $\varphi\in C_0^\infty (Q_\varrho(z_0) )$, affine function $l(z)=l(x)$ such that $|l(z_0)|+|\nabla l|\le M$ we have:
\begin{equation}\label{3.1}
\left|\diagintop_{Q_\varrho(z_0)}(u-l)\varphi_{,t}-{\partial A\over\partial q}(z_0,l(z_0),\D l)
(\D u-\D l)\D\varphi\right|
\le C_{lin} (M) \left[ \omega(M+1,\phi_p) \phi_2^\frac{1}{2} + \phi_p + \psi_p+\varrho^\beta \right] \sup_{Q_\varrho(z_0)}|\D\varphi|
\end{equation}
\end{lem}

\begin{proof}
Use time-independence of $l$ to get from weak formulation of (\ref{1.1}) that for any $\varphi\in C_0^\infty (Q_\varrho(z_0) )$ holds
\begin{equation}\label{3.2}
\begin{aligned}
0&=\diagintop_{Q_\varrho(z_0)}(u-l) \varphi_{,t}-A(z,u,\D u)\D\varphi
\end{aligned}
\end{equation}
which  by adding and subtracting certain terms yields
\begin{equation}\label{3.25}
\begin{aligned}
&\diagintop_{Q_\varrho(z_0)}(u-l)\varphi_{,t}-{\partial A\over\partial q}(z_0,l(z_0),\D l)(\D u-\D l)\D\varphi=\\
&\diagintop_{Q_\varrho(z_0)}[A(z,u,\D u)-A(z_0,l(z_0),\D u)]\D\varphi
+\diagintop_{Q_\varrho(z_0)}\left[A(z_0,l(z_0),\D u)-{\partial A\over\partial q}(z_0,l(z_0),\D l)\cdot(\D u-\D l) \right] \D\varphi.
\end{aligned}
\end{equation}
To obtain our thesis we need to estimate the right-hand-side of (\ref{3.25}).
First, estimate second integral on the r.h.s. of (\ref{3.25}) with respect to 
the splitting of $Q_\varrho(z_0)$ into 
\begin{equation}\label{3.27}
\begin{aligned}
Q_\varrho^s  &= Q_\varrho(z_0) \cap \{|\D u-\D l|\le1\}\\
Q_\varrho^b &= Q_\varrho(z_0)\cap \{|\D u-\D l|>1\}
\end{aligned}
\end{equation}
i.e. 
\begin{equation}\label{3.28}
{1\over|Q_\varrho|}  \int_{Q^s_\varrho}\left[A(z_0,l(z_0),\D u)-{\partial A\over\partial q}(z_0,l(z_0),\D l)\cdot(\D u-\D l) \right] \D\varphi
\end{equation}
and
\begin{equation}\label{3.29}
{1\over|Q_\varrho|}   \int_{Q^b_\varrho}\left[A(z_0,l(z_0),\D u)-{\partial A\over\partial q}(z_0,l(z_0),\D l)\cdot(\D u-\D l) \right] \D\varphi
\end{equation}
without loss of generality assume that neither $Q_\varrho^s$ nor $Q_\varrho^b$ is empty.\\
Since $A(\cdot,0)=\dashint_{Q^s_\varrho}A(z_0,l(z_0),\D l)\D\varphi=0$ one infers that \eqref{3.28} is
\begin{equation}\label{3.3}
{1\over|Q_\varrho|}  \int_{Q^s_\varrho}\intop_0^1\left[ {\partial A\over\partial q}(z_0,l(z_0),\D l+\tau(\D u-\D l))-
{\partial A\over\partial q}(z_0,l(z_0),\D l)\right] (\D u-\D l)\D\varphi d\tau.
\end{equation}
On $Q_\varrho^s$ holds $|l(z_0)|+|\D l|+|\D u-\D l|\le M+1$ in view of assumptions on $l$, so by (\ref{1.5})
\begin{equation}\label{3.32}
\begin{aligned}
&\left|\intop_0^1\left[{\partial A\over\partial q}(z_0,l(z_0),\D l+\tau(\D u-\D l))-
{\partial A\over\partial q}(z_0,l(z_0),\D l)\right](\D u-\D l)\D\varphi d\tau\right|\\
&\le C (M) \;  \omega(M+1,|\D u-\D l|^p)
|\D u-\D l| |\D\varphi|.
\end{aligned}
\end{equation}
Merging \eqref{3.3} and \eqref{3.32} one has
\begin{equation}\label{3.4}
\begin{aligned}
&{1\over|Q_\varrho|} \intop_{Q^s_\varrho}\intop_0^1\left[ {\partial A\over\partial q}(z_0,l(z_0),\D l+\tau(\D u-\D l))-
{\partial A\over\partial q}(z_0,l(z_0),\D l)\right] (\D u-\D l)\D\varphi d\tau\\
&\le C (M) \left[\dashint_{Q_\varrho(z_0)} \;\omega^{p}(M+1, |\D u-\D l|^p)\right]^{1/p}
\left[ \dashint_{Q_\varrho(z_0)}|\D u-\D l|^{p'} \right]^{1/{p'}} \sup_{Q_\varrho(z_0)} |\D\varphi|  \\
& \le C(M) \; \omega (M+1, \phi_p) \phi^{1/{p'}}_{p'}  \sup_{Q_\varrho(z_0)} |\D\varphi|  \le C(M) \; \omega (M+1, \phi_p) \phi^{1/{2}}_{2}  \sup_{Q_\varrho(z_0)} |\D\varphi| 
\end{aligned}
\end{equation}
where the last two inequalities hold by concavity of $\omega^{p}(t,\cdot)$ and $p \ge 2$ ( this is in fact the only place here where we use assumption for $p$). Therefore we can estimate \eqref{3.28} as follows
\begin{equation}\label{qs}
{1\over|Q_\varrho|}  \int_{Q^s_\varrho}\left[A(z_0,l(z_0),\D u)-{\partial A\over\partial q}(z_0,l(z_0),\D l)\cdot(\D u-\D l) \right] \D\varphi  \le C(M) \; \omega (M+1, \phi_p) \phi^{1/p}_{p}  \sup_{Q_\varrho(z_0)} |\D\varphi| 
\end{equation}
Consider now (nonempty) $Q_\varrho^b$. One has for any $s>1$
\begin{equation}\label{3.5}
{|Q_\varrho^b|\over|Q_\varrho|} \le 1 \le  \diagintop_{Q_\varrho^b}|\D u-\D l|^s
\end{equation}
because 
\[
|Q_\varrho^b|\le\intop_{Q_\varrho^b}|\D u-\D l|\le\left(\intop_{Q_\varrho^b}|\D u-\D l|^s\right)^{1/s}
|Q_\varrho^b|^{1/s'}\le\left(\dashint_{Q_\varrho^b}|\D u-\D l|^s\right)^{1/s}
|Q_\varrho^b|^{1/s'}|Q_\varrho|^{1/s}
\]
From $|\D l|+|l(z_0)|\le M$ and (\ref{1.6}), (\ref{1.7}) we estimate  \eqref{3.29}
\begin{equation}\label{3.6}
\begin{aligned}
&{1\over|Q_\varrho|} \intop_{Q_\varrho^b}\left[A(z_0,l(z_0),\D u)-{\partial A\over\partial q}
(z_0,l(z_0),\D l)(\D u-\D l)\right]\D\varphi \le\\
&\sup_{Q_\varrho(z_0)}|\D\varphi|{ C(M) \over|Q_\varrho|}\intop_{Q_\varrho^b} 1+|\D u|^{p-1}+ |\D u-\D l|\le
\sup_{Q_\varrho(z_0)} |\D\varphi|{ C(M) \over|Q_\varrho|}\intop_{Q_\varrho^b}1+|\D u-\D l|^{p-1} + |\D u-\D l| \\
&\le C(M) \sup_{Q_\varrho(z_0)}|\D\varphi|\left[{ |Q_\varrho^b|\over|Q_\varrho|}+\left(\diagintop_{Q_\varrho}
|\D u-\D l|^p\right)^{1/p'} {|Q_\varrho^b|^{1/p} \over|Q_\varrho|^{1/p}} + \left(\diagintop_{Q_\varrho}
|\D u-\D l|^p\right)^{1/p} {|Q_\varrho^b|^{1/p'} \over|Q_\varrho|^{1/p'}}\right] \\
&\le C(M) \sup_{Q_\varrho(z_0)}|\D\varphi| \phi_p
\end{aligned}
\end{equation}
where the last inequality holds in view of (\ref{3.5}) with $s=p$ and  $s=p'$. Combine estimates \eqref{qs} and \eqref{3.6} to get
\begin{equation}\label{linA}
\diagintop_{Q_\varrho(z_0)}\left[A(z_0,l(z_0),\D u)-{\partial A\over\partial q}(z_0,l(z_0),\D l)\cdot(\D u-\D l) \right] \D\varphi \le C(M) [\omega (M+1, \phi_p) \phi^{1/p}_{p} + \phi_p ] \sup_{Q_\varrho(z_0)} |\D\varphi|.
\end{equation}
It remains to estimate the first term in (\ref{3.25}); use \eqref{alg_1} with $P \equiv \D u$ to get (\ref{1.8}) 
\begin{equation}\label{3.7}
|A(z,u,\D u)-A(z_0,l(z_0),\D u)| \le  C(M) \varrho^\beta
\left [1+  |\D u-\D l|^{p} + \left| \frac{u - l}{\varrho} \right|^p \right]
\end{equation}
Inequality \eqref{3.7} used to estimate the first term of the right-hand-side of (\ref{3.2}) gives
\begin{equation}\label{linB}
\begin{aligned}
&\left|\diagintop_{Q_\varrho(z_0)}(A(z,u,\D u)-A(z_0,l(z_0),\D u))\D\varphi\right| \le \sup_{Q_\varrho(z_0)}|\D\varphi|C(M) \varrho^\beta
\left [1+ \psi_p + \phi_p \right]
\end{aligned}
\end{equation}
Inequalities (\ref{linA}), (\ref{linB}) used in (\ref{3.25}) give thesis.
\end{proof}
\begin{lem}\label{lem4.1}
(Local inequalities) Take $p \ge 2$. Let $u \in C(-T, 0; L^2(\Omega)) \cap L^p(-T,0; W^{1,p}(\Omega))$ be a weak solution to (\ref{1.1}) with structure conditions (\ref{1.6} --- \ref{1.2b}). Then the following inequalities hold for any $Q_\varrho(z_0) \subset Q$ with $\varrho \le 1$ and constants being nondecreasing functions of their parameters
\begin{multline}\label{pre4.1.1}
\left| \intop_{B_\varrho(x_0)} (u (t,x) - u(\tau,x) )\eta_\varrho(x) dx \right| \le \\
C_{\eqref{pre4.1.1} } ( |(u)_{z_0}| + |(\D u)_{z_0}|) \varrho  \diagintop_{Q_\varrho(z_0)}  \left( \varrho^\beta + |u- (u)_{z_0} |^\beta+ |u- (u)_{z_0} | + |\D u-\ (\D u)_{z_0} |^{p-1} +  |\D u-\ (\D u)_{z_0} | \right),
\end{multline}
\begin{equation}\label{prepre4.1.1}
\left| \intop_{B_\varrho(x_0)} (u (t,x) - u(\tau,x) )\eta_\varrho(x) dx \right| \le  C_{\eqref{prepre4.1.1} } (|(u)_{z_0}| + |(\D u)_{z_0}|)  \varrho  \diagintop_{Q_\varrho(z_0)}  \left(1+ |\D u|^{p-1} \right) \end{equation}
for $t, \tau \in (t_0 - \varrho^2, t_0)$, where $\eta_\varrho (x)$ denotes a standard mollifier in space, supported in $B_\varrho(x_0)$;
\begin{equation}\label{4.1.1}
\sup_{t \in \left( t_0 - \left(\frac{\varrho}{2} \right)^2, t_0 \right)} |u-l|^2_{L^2 (B_\frac{\varrho}{2} )} +  \diagintop_{Q_{\varrho\over2}(z_0)}|\D u-\D l|^2 +|\D u-\D l|^p \le C_{Cacc} (M) \left[\diagintop_{Q_\varrho(z_0)}
\left|{u-l\over\varrho}\right|^2 + \diagintop_{Q_\varrho(z_0)}
\left|{u-l\over\varrho}\right|^p + \varrho^{2\beta}\right],
\end{equation}
\begin{equation}\label{4.1.1'}
\tag{\ref{4.1.1}$'$}
\diagintop_{Q_{\varrho\over2}(z_0)}|\nabla u-\nabla l|^2 \le C'_{Cacc} (M)\left[\diagintop_{Q_\varrho(z_0)}
\left|{u-l\over\varrho}\right|^2 + \diagintop_{Q_\varrho(z_0)}
\left|{u-l\over\varrho}\right|^p + \varrho^{2\beta}\right]
\end{equation}
where $l$ is an affine function depending only on $x$ and satisfying $|l(z_0)|+|\nabla l|\le M$ and $\beta \in (0,1)$ is given by (\ref{1.8}).
\end{lem}

\begin{proof}
Fix arbitrary numbers $t, \tau$ and nonnegative $ \varepsilon, \tilde \varepsilon$ satisfying
\begin{equation}
t_0 - \varrho^2
\le t <  t+\tilde \varepsilon < \tau - \varepsilon < \tau \le t_0 \end{equation}
and the continuous, piecewise affine cutoff function $\sigma(s)  \in [0,1]$ defined by
\begin{equation}\label{sigma}
\begin{aligned}
\sigma_{t, \tau, \varepsilon, \tilde \varepsilon} (s) &=\begin{cases}1&  \text{on} \; (t + \tilde \varepsilon, \; \tau-\varepsilon ),\\
0 &  \text{on} \; (t, \; \tau)^c,\end{cases}\\
\sigma_{t, \tau, \varepsilon, \tilde \varepsilon}'(s) &=\begin{cases} 1/  \tilde \varepsilon &  \text{on} \;  (t, \; t + \tilde \varepsilon ),\\
- 1/ \varepsilon &  \text{on} \;  (\tau -  \varepsilon, \; \tau ),\\
0 &  \text{otherwise}.\end{cases}\\
\end{aligned}
\end{equation}
Let us first show \eqref{pre4.1.1}.  Test (\ref{1.1}) with $\sigma_{t, \tau, \varepsilon, \varepsilon} \eta_\varrho $, obtaining
\begin{equation}\label{mean_ac}
\intop_{B_\varrho(x_0)} (u (t,x) - u(\tau,x) )\eta_\varrho(x) = \int^t_\tau \intop_{B_\varrho(x_0)} A(z,u,\D u) \D \eta_\varrho (x) ds
\end{equation}
by sending $\varepsilon \to 0$ (this holds pointwisely in time, because $u \in C(L^2)$). Estimate the r.h.s. of \eqref{mean_ac} using that  $|\nabla \eta_\varrho(x)| \le C \varrho^{-(d+1)} $
\begin{multline}\label{lok_rz}
 \left| \int^t_\tau \intop_{B_\varrho(x_0)} A(z,u,\D u) \D \eta_\varrho(x) dx ds \right| =  \left| \int^t_\tau \intop_{B_\varrho(x_0)} (A(z,u,\D u) - A(z_0,(u)_{z_0} , (\D u)_{z_0} )) \D \eta_\varrho(x) dx ds \right| \\
 \le C \intop_{Q_\varrho(z_0)}  \left| A(z,u,\D u) - A(z_0,(u)_{z_0} , (\D u)_{z_0} )\right| |\varrho|^{-(d+1)}  \le C \varrho \diagintop_{Q_\varrho(z_0)}  \left| A(z,u,\D u) \pm A(z,u, (\D u)_{z_0} )  - A(z_0 ,(u)_{z_0} , (\D u)_{z_0} )\right| \\  \le C ( |(u)_{z_0}| + |(\D u)_{z_0}|) \varrho  \diagintop_{Q_\varrho(x_0)}  \left( \varrho^\beta + |u- (u)_{z_0} |^\beta+ |u- (u)_{z_0} | + |\D u-\ (\D u)_{z_0} |^{p-1} +  |\D u-\ (\D u)_{z_0} | \right) 
\end{multline}
where the last inequality comes from adding estimates \eqref{alg_3}, \eqref{alg_2} with $P \equiv \D u$ and $ l (x) \equiv (u)_{z_0} + (\D u)_{z_0} (x-x_0)$. This ends the proof of  \eqref{pre4.1.1}.  To get  \eqref{prepre4.1.1}, when estimating \eqref{lok_rz},  we use inequality \eqref{alg_4} instead of \eqref{alg_3} and \eqref{alg_2}.   \\
Let us now turn our attention to the energy estimate  \eqref{4.1.1}. To show it, choose a smooth cutoff function $\theta(x) \in [0,1]$ satisfying
\begin{equation}\label{theta} 
\begin{aligned}
\theta(x) &=\begin{cases}1& \text{on} \; B_{\varrho/2}(z_0),\\ 
0&  \text{on} \;  B_\varrho^c(z_0),\end{cases}\\
|\nabla\theta| &\le{4/ \varrho}.
\end{aligned}
\end{equation}
Test (\ref{1.1}) with $\varphi= (u-l) \theta^2(x) \sigma_{t_0 -\varrho^2, \tau, \varepsilon, \frac{3}{4} \varrho^2} (s) $, suppressing for now parameters of cutoff function in time, thus writing $\sigma$. The evolutionary part yields
\begin{equation}\label{4.1.2}
\begin{aligned}
&\intop_{Q_\varrho(z_0)}u\varphi_{,s} = \intop_{Q_\varrho(z_0)}(u- l )\varphi_{,s} =\intop_{Q_\varrho(z_0)}(u-l)\theta^2(\sigma_{,s}(u-l)+\sigma(u-l)_{,s})=
\intop_{Q_\varrho(z_0)} |u-l|^2\theta^2\sigma_{,s}\\
&\quad+ {1\over2} \intop_{Q_\varrho(z_0)} (|u-l|^2 \theta^2)_{,s}\sigma= {1\over2} \intop_{Q_\varrho(z_0)}
|u-l|^2\theta^2\sigma_{,s}\le \intop_{B_\varrho\times\left(t_0-\varrho^2,t_0-{\varrho^2\over 4}\right)}
\left|{u-l\over\varrho}\right|^2 \theta^2 -\frac{1}{\varepsilon} \intop_{B_\varrho\times\left(\tau-\varepsilon, \tau \right)}
\left|{u-l}\right|^2 \theta^2
\end{aligned}
\end{equation}
the last inequality holds, because $|\sigma' |  \le 2 \varrho^{-2} $ on $ \left(t_0-\varrho^2,t_0- \varrho^2 / 4 \right)$ in view of \eqref{sigma}.
As $\intop_{Q_\varrho(z_0)}A(z_0,l(z_0),\D l)\D \varphi =0$, for the main part holds
\begin{equation}\label{4.1.5}
\begin{aligned}
&\intop_{Q_\varrho(z_0)}A(z,u,\D u)\D \varphi =\intop_{Q_\varrho(z_0)}[A(z,u,\D u)- A(z,u,\D l)]\D\varphi
+ [A(z,u,\D l)-A(z_0,l(z_0), \D l)]\D\varphi
=\\
& \intop_{Q_\varrho(z_0)}[A(z,u,\D u)-A(z,u,\D l)]
\D(u-l)\theta^2\sigma + [A(z,u,\D u)-A(z,u,\D l)]2\theta \nabla\theta \colon (u-l)\sigma 
+ [A(z,u,\D l) - A(z_0,l(z_0),\D l)]\D\varphi.
\end{aligned}
\end{equation}
For $p \ge 2$ holds
\begin{equation}\label{4.1.51}
|A-B|^2 \int_0^1 ( 1 + | s (A - B)  +  B |^2)^\frac{p-2}{2}  ds \ge c (|A-B|^2 +|A-B|^p)
\end{equation}
hence, using the  assumption \eqref{1.2b} one obtains
\begin{multline} \label{4.1.55}
[A(z,u,\D u)-A(z,u,\D l)] \D(u-l) = \int_0^1 \frac{\partial A}{\partial q}  (z, u, s ( \D u - \D l ) + \D l) ( \D u - \D l ) \cdot  ( \D u - \D l ) ds \\
\ge \int_0^1  \lambda (1 + | s ( \D u - \D l ) + \D l |^2)^{\frac{p-2}{2}} |  \D u - \D l  |^2 ds \ge c ( |  \D u - \D l  |^p + |  \D u - \D l  |^2 )
\end{multline}
\noindent
Inequalities \eqref{4.1.2},   \eqref{4.1.5}, \eqref{4.1.55} show that testing (\ref{1.1}) with $\varphi=\theta^2\sigma (u-l)$ yields the following estimate
\begin{multline}\label{4.1.6}
\frac{1}{\varepsilon}  \intop_{B_\varrho\times\left(\tau-\varepsilon, \tau \right)} \left|{u-l}\right|^2 \theta^2
 + \intop_{Q_\varrho(z_0)}|\D(u-l)|^2\theta^2\sigma + \intop_{Q_\varrho(z_0)}|\D(u-l)|^p\theta^2\sigma  \le C \intop_{Q_\varrho(z_0)}
\left|{u-l\over\varrho}\right|^2 +\\
C \intop_{Q_\varrho(z_0)}  \left|  [A(z,u,\D u)-A(z,u,\D l)]2\theta \nabla\theta \colon (u-l)\sigma \right| 
+C \intop_{Q_\varrho(z_0)} \left|  [A(z,u,\D l)-A(z_0,l(z_0),\D l)]\D\varphi \right| \equiv  C \intop_{Q_\varrho(z_0)}
\left|{u-l\over\varrho}\right|^2 + I + II 
\end{multline}
Let us estimate $I$ by (\ref{alg_2}) with $P \equiv \D u$ getting for $\varrho \le 1$
\begin{multline}
I \le C(M) \intop_{Q_\varrho(z_0)}\left( \left|u-l\right|^\beta+ |u-l| + |\D u-\D l|^{p-1}  \right) \left| \frac{u-l}{\varrho} \right| \sigma \theta 
 \le \\
  \frac{1}{2} \intop_{Q_\varrho(z_0)} |\D u-\D l|^p \theta^2\sigma + C(M)
  \intop_{Q_\varrho(z_0)}
\left|{u-l\over\varrho}\right|^2 + \left|{u-l\over\varrho}\right|^p + {|u-l|^{\beta+1}\over\varrho}.
\end{multline}
Use \eqref{alg_3} to obtain
\begin{equation}
 |A(z_0,l(z_0), \D l) - A(z, u, \D l)| |\D\varphi | \le C(M) \varrho^\beta \left[1+ |u-l|^\beta \right] \left({|u-l|\over\varrho}+|\D u-\D l|\right)  \theta \sigma
 \end{equation}
with which we estimate $II$ 
\begin{multline}
II \le C(M)\intop_{Q_\varrho(z_0)}\left[\varrho^\beta|\D u-\D l|+\varrho^\beta{|u-l|\over\varrho}+
{|u-l|^{\beta+1}\over\varrho}+|u-l|^\beta|\D u-\D l|\right]  \theta \sigma\\
\le \frac{1}{2} \intop_{Q_\varrho(z_0)}|\D u-\D l|^2  \theta^2 \sigma +C(M) \intop_{Q_\varrho(z_0)}
\left[\varrho^{2\beta}+|u-l|^{2\beta}+\left|{u-l\over\varrho}\right|^2+{|u-l|^{\beta+1}\over\varrho}\right]
\end{multline}
Estimates for $I$ and $II$ give together
\begin{equation}\label{I+II}
I+II \le   \frac{1}{2} \intop_{Q_\varrho(z_0)} \left[ |\D u-\D l|^p +  |\D u-\D l|^2 \right] \theta^2\sigma + C(M) \intop_{Q_\varrho(z_0)}
\left[\varrho^{2\beta}+|u-l|^{2\beta}+\left|{u-l\over\varrho}\right|^2 +\left|{u-l\over\varrho}\right|^p+{|u-l|^{\beta+1}\over\varrho}\right]
\end{equation}
In view of  $\beta<1, \varrho \le 1$ one has
\begin{equation*}
\intop_{Q_\varrho(z_0)}|u-l|^{2\beta}\le\varrho^{2\beta}\intop_{Q_\varrho(z_0)}
\left|{u-l\over\varrho }\right|^{2\beta}\le  \varrho^{2 \beta} |Q_\varrho|+ C\intop_{Q_\varrho(z_0)}
\left|{u-l\over\varrho}\right|^2,
\end{equation*}
\begin{equation*}
\intop_{Q_\varrho(z_0)}{|u-l|^{\beta+1}\over\varrho}\le\varrho^{{2\beta} \over {1-\beta}}|Q_\varrho|+\intop_{Q_\varrho(z_0)}
\left|{u-l\over\varrho}\right|^2 \le \varrho^{{2\beta}}|Q_\varrho|+\intop_{Q_\varrho(z_0)}
\left|{u-l\over\varrho}\right|^2
\end{equation*}
Consequently, \eqref{I+II} takes the form
\begin{equation}\label{I+II'}
I+II \le   \frac{1}{2} \intop_{Q_\varrho(z_0)} \left[ |\D u-\D l|^p +  |\D u-\D l|^2 \right] \theta^2\sigma + C(M) \left[\varrho^{2\beta} |Q_\varrho|+ \intop_{Q_\varrho(z_0)}
\left|{u-l\over\varrho}\right|^2 +\left|{u-l\over\varrho}\right|^p \right].
\end{equation}
hence (\ref{4.1.6}) with \eqref{I+II'} yields
\begin{multline}\label{cacc_mezz}
\frac{1}{\varepsilon}  \intop_{B_{\varrho/2}(z_0)\times\left(\tau-\varepsilon, \tau \right)} \left|{u-l}\right|^2 
 + \intop_{t_0-\left({\varrho\over2}\right)^2}^{\tau-\varepsilon}\intop_{B_{\varrho/2}(z_0)}
 |\D(u-l)|^p + |\D(u-l)|^2 \le
C(M) \left[\varrho^{2\beta}
|Q_\varrho|+\intop_{Q_\varrho(z_0)}\left|{u-l\over\varrho}\right|^2 + \left|{u-l\over\varrho}\right|^p \right].
\end{multline}
First, use inequality \eqref{cacc_mezz} for $\tau = t_0$, neglecting the first term of the left-hand-side. This estimate is uniform in $\varepsilon$, so we obtain
 \begin{equation}\label{4.1.12}
 \dashint_{Q_{\varrho\over2}(z_0)}|\D u-\D l|^2 + \dashint_{Q_{\varrho\over2}(z_0)}|\D u-\D l|^p \le C(M)\left[\diagintop_{Q_\varrho(z_0)}
\left|{u-l\over\varrho}\right|^2 + \diagintop_{Q_\varrho(z_0)}
\left|{u-l\over\varrho}\right|^p + \varrho^{2\beta}\right]
\end{equation}
Next, drop second part of  left-hand-side of \eqref{cacc_mezz} and consider any $\tau $ in interval of admissibility $\left(t_0 - \frac{\rho^2}{4}, t_0 \right]$; this via Steklov averages argument gives rise to
\begin{equation}\label{4.1.13}
\sup_{t \in \left( t_0 - \left(\frac{\varrho}{2} \right)^2, t_0 \right)} |u-l|^2_{L^2 (B_\frac{\varrho}{2} )}  \le C (M)\left[\diagintop_{Q_\varrho(z_0)}
\left|{u-l\over\varrho}\right|^2 + \diagintop_{Q_\varrho(z_0)}
\left|{u-l\over\varrho}\right|^p + \varrho^{2\beta}\right].
\end{equation}
Combining \eqref{4.1.12} and \eqref{4.1.13} we have the first Caccioppoli estimate  \eqref{4.1.1}. It implies, in conjunction with the Korn's inequality \eqref{1.3.1} used for $(u-l) (t)$, the following estimate 
\begin{multline}\label{4.1.1'pf}
\diagintop^{t_0}_{t_0 - {\varrho^2\over4}} \diagintop_{B_{\varrho\over2}(x_0)}|\nabla u-\nabla l|^2 \le  \diagintop^{t_0}_{t_0 - {\varrho^2\over4}}  K_2  \left[\diagintop_{B_{\varrho\over2}(x_0)}|\D u-\D l|^2 + \diagintop_{B_{\varrho\over2}(x_0)}
\left|{u-l\over{\varrho/2}}\right|^2  \right] \le \\
 K_2  C_{Cacc} (M)\left[\diagintop_{Q_\varrho(z_0)}
\left|{u-l\over\varrho}\right|^2 + \diagintop_{Q_\varrho(z_0)}
\left|{u-l\over\varrho}\right|^p + \varrho^{2\beta}\right] + K_2 2^{n+4} \diagintop_{Q_{\varrho}(z_0)}
\left|{u-l\over\varrho}\right|^2, 
\end{multline}
which justifies \eqref{4.1.1'} with $C'_{Cacc} \equiv K_2 (C_{Cacc} (M) +2^{n+4})$
\end{proof}
Next, we restate the linearization lemma (Lemma \ref{lem3.1}) using local inequalities of Lemma \ref{lem4.1} in a way useful for further computations. To proceed, introduce the following useful quantities
\begin{mydef}\label{energies}
$ E_{z_0, l} (\varrho)$ denotes the $L^2-L^p$ excess energy 
\begin{equation}\label{4.15.1'}
 E_{z_0, l} (\varrho) \equiv \psi_{2, z_0, l}   (\varrho) + \psi_{p , z_0, l} (\varrho) \qquad \left(= \diagintop_{Q_{z_0}(\varrho)} \left|{u-l\over\varrho}\right|^2 + \diagintop_{Q_{z_0}(\varrho)} \left|{u-l\over\varrho}\right|^p  \right)
 \end{equation}
and $ \tilde  E_{z_0, l} (\varrho)$ denotes the perturbed $L^2-L^p$ excess energy
\begin{equation}\label{4.15.1''}
\tilde  E_{z_0, l} (\varrho) \equiv  E_{z_0, l} (\varrho) + \varrho^{2 \beta} \quad 
\end{equation}
where $\psi_{2, z_0, l}   (\varrho), \; \psi_{p , z_0, l}   (\varrho)$ are given as in Definition \ref{inxs}.
 \end{mydef}
\begin{mydef}\label{vau}
Introduce normalization factor $\gamma$, which depends on parameters $ \delta, \varrho, l$
\begin{equation}\label{4.15.2'}
\gamma_{l,\delta} (\varrho) \equiv \sqrt{ E_{z_0, l} (\varrho) + (\delta/2)^{-2}\varrho^{2\beta}} 
 \end{equation}
\end{mydef}

\begin{cor}\label{cor4.2}
Take $p \ge 2$ and fix $M$. Let $u \in C(-T, 0; L^2(\Omega)) \cap L^p(-T,0; W^{1,p}(\Omega))$ be a weak solution to (\ref{1.1}) with structure conditions (\ref{1.2a} --- \ref{1.2b}), (\ref{1.7} --- \ref{1.5}). There exists such constant $C_{\ref{cor4.2}}(M)$ that for any affine function $l$  depending only on $x$ and satisfying $|l(x_0)|+|\nabla l|\le M$ and any $\delta \in (0,1)$ hold for 
\begin{equation}\label{4.15.2''}
v \equiv \frac{u-l}{C_{\ref{cor4.2}}(M)  \gamma_{l,\delta} (\varrho)}
\end{equation}
the following inequalities
\begin{equation}\label{4.2.1}
\left|\diagintop_{Q_{\varrho/2}(z_0)}v\varphi_{,t}-{\partial A\over\partial q}(z_0,l(z_0),\D l)\D v\D\varphi\right|
\le  \left[ \omega \left(M+1, \tilde  E_{z_0, l}  (\varrho) \right) +   \tilde E_{z_0, l} ^\frac{1}{2}  (\varrho) + \delta/2 \right] 
\sup_{Q_{\varrho/2}(z_0)}|\D\varphi|
\end{equation}
and
\begin{equation}\label{4.2.2}
\dashint_{Q_{\varrho/2}(z_0)} \left|\frac{v}{\varrho / 2} \right|^2 +  |\nabla v|^2+( {C_{\ref{cor4.2}}(M)  \gamma_{l,\delta} (\varrho)})^{p-2} \left[ \left|\frac{v}{\varrho / 2} \right|^p +  |\nabla v|^p \right] \le 1,
\end{equation}
where $Q_\varrho(z_0) \subset Q$ is an arbitrary local cylinder  with $\varrho \le 1$.
\end{cor}
\begin{proof}
We suppress parameters of the excess energies writing $E, \tilde E$ for $ E_{z_0, l},  \tilde E_{z_0, l}$ and similarily for moments $\psi, \phi$. Take
\begin{equation}\label{c_vau}
C_{\ref{cor4.2}}(M) \equiv  2^{n+2+p}  (1 + \max(K_2, K_p))^{1/2} C_{lin} (M) \left(1+C^{1+\frac{1}{p}}_{Cacc} (M) \right)
\end{equation}
where $K_2, K_p$ are constants from the Korn's inequality \eqref{1.3.1} of Lemma \ref{lem1.3}; with this choice of $C_{\ref{cor4.2}}(M),  $ \eqref{4.2.3} yields \eqref{4.2.1}. \\
Assumptions of Lemmas \ref{lem3.1}, \ref{lem4.1} are fulfilled.
Linearization inequality
(\ref{3.1}) with a Caccioppoli estimate (\ref{4.1.1}) give
\begin{multline}\label{4.2.3}
\left|\diagintop_{Q_{\varrho/2}(z_0)}(u-l)\varphi_{,t}-{\partial A\over\partial q}(z_0,l(z_0),\D l)
(\D u-\D l)\D\varphi\right|
\le \\
C_{lin} (M) \left[ \omega \left(M+1,\phi_p  \left( \frac{\varrho}{2} \right) \right) \phi_2^\frac{1}{2}  \left( \frac{\varrho}{2} \right) + \phi_p  \left( \frac{\varrho}{2} \right) + \psi_p  \left( \frac{\varrho}{2} \right) +  \left( \frac{\varrho}{2} \right)^\beta \right] \sup_{Q_{\varrho/2}(z_0)}|\D\varphi| \\ 
\le
C_{lin} (M) C_{Cacc} (M) \left[ \omega \left(M+1, C_{Cacc} (M)   \tilde E (\varrho) \right)  \tilde E^\frac{1}{2}  (\varrho) +  \tilde E (\varrho)  +  2^{n+2+p} E (\varrho)  + \varrho^\beta \right] \sup_{Q_{\varrho/2}(z_0)}|\D\varphi| \le \\ 
 2^{n+2+p}  C_{lin} (M) C^{1+\frac{1}{p}}_{Cacc} (M) \gamma_{l,\delta} (\varrho) \left[ \omega \left(M+1, \tilde  E (\varrho) \right) +   \tilde E^\frac{1}{2}  (\varrho) + \delta/2 \right] \sup_{Q_{\varrho/2}(z_0)}|\D\varphi|
\end{multline}
The last inequality holds by concavity of $\omega^p$ with respect to its second variable giving for $c>1$:  $\omega(M+1,c\alpha)\le c^\frac{1}{p}\omega(M+1,\alpha)$, and by definition \eqref{4.15.2'} of $\gamma_{l,\delta} (\varrho)$. Consequently we have (\ref{4.2.1}).
 Let us now justify inequality (\ref{4.2.2}). Using the Korn's inequality  \eqref{1.3.1} from Lemma \ref{lem1.3} compute
\begin{multline}\label{do6.35}
\dashint_{Q_{\varrho/2}(z_0)} \left|\frac{v}{\varrho / 2} \right|^2 +  |\nabla v|^2+( {C_{\ref{cor4.2}}(M)  \gamma_{l,\delta} (\varrho)})^{p-2} \left[ \left|\frac{v}{\varrho / 2} \right|^p +  |\nabla v|^p \right] \le \\
(1 + \max(K_2, K_p)) \dashint_{Q_{\varrho/2}(z_0)} \left|\frac{v}{\varrho / 2} \right|^2 +  |\D v|^2+( {C_{\ref{cor4.2}}(M)  \gamma_{l,\delta} (\varrho)})^{p-2} \left[ \left|\frac{v}{\varrho / 2} \right|^p +  |\D v|^p \right] \\ \le
\left(  2^{n+2+p}  C_{lin} (M) \left(1+C^{1+\frac{1}{p}}_{Cacc} (M) \right)  \gamma_{l,\delta} (\varrho) \right)^{-2} \dashint_{Q_{\varrho/2}(z_0)}|\D(u-l)|^2+ \left|{u-l\over{\varrho/2}}\right|^2 + |\D(u-l)|^p+ \left|{u-l\over{\varrho/2}}\right|^p \\ \le
\left( 2^{n+2+p}  C_{lin} (M) \left(1+C^{1+\frac{1}{p}}_{Cacc} (M) \right) \right)^{-2} {\tilde E_{z_0, l}}^{-1} (\varrho) \dashint_{Q_{\varrho/2}(z_0)}|\D(u-l)|^2+ \left|{u-l\over{\varrho/2}}\right|^2 + |\D(u-l)|^p+ \left|{u-l\over{\varrho/2}}\right|^p  \\ \le
\left( 2^{n+2+p}  C_{lin} (M) \left(1+C^{1+\frac{1}{p}}_{Cacc} (M) \right) \right)^{-2} ( 2^{n+2+p} + C_{Cacc} (M))
\end{multline}
The last three inequalities come, respectively, from: definition \eqref{4.15.2''} of $v$ and choice \eqref{c_vau} of constant  $C_{\ref{cor4.2}}(M) $;  the Definition \ref{vau} of $ \gamma_{l, \delta} (\varrho)$ and the fact that $\delta \le 1$; the Caccioppoli inequality \eqref{4.1.1} and the Definition \ref{energies} of perturbed excess energy $\tilde E$. As $C_{lin}, C_{Cacc}$ are bigger than $1$, \eqref{do6.35} implies (\ref{4.2.2}).
\end{proof}

\section{Partial regularity}\label{partial}
\setcounter{equation}{0}

\noindent
First we merge the local inequality of Corollary \ref{cor4.2} and the caloric approximation into a building block of a further partial regularity result. Recall from Lemma \ref{lem1.2} that  $l^{(s)}_{z_0, \varrho} (x)$ is the affine function, depending only on space variable, which minimizes $\dashint_{Q_\varrho(z_0)}|u-l|^s$; $l^{(s)}_{\varrho} (x)$ denotes this function, when dependence on $z_0$ is irrelevant.
\begin{lem}\label{lem6.1}
Let $p \ge 2$ and $u \in C(-T, 0; L^2(\Omega)) \cap L^p(-T,0; W^{1,p}(\Omega))$  be a weak solution to (\ref{1.1}) under structure conditions (\ref{1.2a} --- \ref{1.5}). Fix  constants $M > 0, \; \alpha \in (0,1) $. Then there exist $\sigma \in (0, 1/4), \; \delta \in (0, 1)$,  such that for any $\varrho < 1$, $z_0 \in Q$ (such that $Q_\varrho (z_0) \subset Q$), we have the following implication.\\ If
\begin{equation}\label{6.1.1}
|l^{(2)}_{z_0, \varrho} (x_0)|+|\nabla l^{(2)}_{z_0, \varrho}  |\le M.
\end{equation}
and
\begin{equation}\label{6.1.2}
\omega(M+1, \tilde  E_{z_0, l^{(2)}_{z_0, \varrho} } (\varrho) )+ \tilde  E^{\frac{1}{2}}_{z_0, l^{(2)}_{z_0, \varrho} } (\varrho)  \le{\delta\over 2 C_{\ref{cor4.2}}(M) },
\end{equation}
then
\begin{equation}\label{shrink_step}
 \tilde  E_{z_0, l^{(2)}_{z_0,{\sigma\varrho}}} \left(\sigma\varrho \right) \le  \sigma^{2\alpha} \left[ E_{z_0, l^{(2)}_{z_0, \varrho} } (\varrho) + \delta^{-2} \varrho^{2\beta} \right]+   ({\sigma \varrho})^{2 \beta}
\end{equation}
\end{lem}

\begin{proof} 

We need certain care to avoid a logical loop. Therefore let us first explicitly define constants:
\begin{equation}\label{6.1.05const}
\begin{aligned}
C_\eqref{2.4.1} (\lambda,  \Lambda, s) & \; \text{ is the constant from \eqref{2.4.1} of Lemma \ref{LinSys} with parameters } (\lambda, \Lambda, K_s, 1, d, s, \infty), \\
C_\eqref{2.4.12} (\lambda,  \Lambda, s) & \; \text{ is the constant from \eqref{2.4.12} of Lemma \ref{LinSys} with parameters } (\lambda, \Lambda, K_s, 0, d, s, s), \\
C_\eqref{6.1.05const} & \equiv C^2_{\ref{cor4.2}}(M) 2^{5p -3} \max_{s \in \{2;p\}} \left( 1+ C_\eqref{2.4.1} (\lambda,  \Lambda, s) + C_\eqref{2.4.12} (\lambda,  \Lambda, s)\right)
\end{aligned}
\end{equation}
we have already fixed in the statement of lemma $M > 0, \; \alpha \in (0,1) $. Now let us fix certain parameters:
\begin{equation}\label{6.1.05}
\begin{aligned}
\sigma &< 1/4\; \text{ so that } \;C_\eqref{6.1.05const} 2^7 \sigma^2 \le  \sigma^{2 \alpha}  \quad  \text{ (which is possible as we have assumed that } \alpha \in (0,1) )\\
\varepsilon & \equiv { (4 \sigma)^{p+ d+2}} 16 \sigma^2,\\
\theta &\equiv 4 \sigma, \\
\Lambda &= C_\eqref{1.7} (M).
\end{aligned}
\end{equation}
Observe that by assumptions \eqref{1.2b}, \eqref{1.4} holds 
\begin{equation}\label{6.1.07}
{\partial A\over\partial q} \left(z_0,l^{(2)}_\varrho(z_0),\D l^{(2)}_\varrho \right) \in S \left(\lambda, \left|{\partial A\over\partial q}(z_0,l^{(2)}_\varrho(z_0),\D l^{(2)}_\varrho) \right| \right) \subset  S \left(\lambda,  \Lambda  \right) 
\end{equation}
i.e. the constant coefficients matrix,  resulting from linearization around $z_0$, belongs to the set of elliptic bilinear, symmetrizing forms as defined in Definition \ref{def5.05}. The imbedding  results from \eqref{1.7} with \eqref{6.1.1}; $\lambda $ is given by \eqref{1.2b} and  $\Lambda $ -- by \eqref{6.1.05}. Consequently, let us fix via Lemma \ref{lem5.1}
\begin{equation}
\delta_\varepsilon
 \equiv  \delta \left( \varepsilon, \lambda, \Lambda  \right) 
 \end{equation}
Take 
\begin{equation}\label{choice_gamma}
\gamma \equiv {C_{\ref{cor4.2}}(M)  \gamma_{l^{(2)}_\varrho, \delta_\varepsilon} (\varrho)}.
\end{equation}
Observe that assumptions of Corollary \ref{cor4.2} are fulfilled; this and assumption (\ref{6.1.2}) give for\begin{equation}\label{def_v}
v \equiv \frac{u-l^{(2)}_\varrho }{C_{\ref{cor4.2}}(M)  \gamma_{l^{(2)}_\varrho, \delta_\varepsilon} (\varrho)} \equiv \frac{u-l}{\gamma},
\end{equation} defined as in \eqref{4.15.2''}, inequalities
\begin{equation}\label{6.1.3}
\left|\diagintop_{Q_{\varrho/2}(z_0)}v\varphi_{,t}-{\partial A\over\partial q}(z_0,l^{(2)}_\varrho(z_0),\D l^{(2)}_\varrho)\D v\D\varphi\right|\le\delta_\varepsilon\sup_{Q_{\varrho/2}(z_0)}|\D\varphi|
\end{equation}
\begin{equation}\label{6.1.45}
\dashint_{Q_{\varrho/2}(z_0)} \left|\frac{v}{\varrho / 2} \right|^2 + \gamma^{p-2} \left|\frac{v}{\varrho / 2} \right|^p+  |\nabla v|^2 + \gamma^{p-2} |\nabla v|^p \le 1
\end{equation}
 By definition \eqref{4.15.2'} and  (\ref{6.1.2}) one has also
\begin{equation} \label{6.1.42}
0 \le \gamma= C_{\ref{cor4.2}}(M)   \sqrt{ E_{z_0, l^{(2)}_\varrho } (\varrho) + (\delta_\varepsilon/2)^{-2}\varrho^{2\beta}}  \le \frac{ C_{\ref{cor4.2}}(M) }{ \delta}  \sqrt{ \tilde  E_{z_0, l^{(2)}_\varrho} (\varrho)  } \le 1.
 \end{equation}
Obervation \eqref{6.1.07} with inequalities (\ref{6.1.3}), (\ref{6.1.45}), \eqref{6.1.42} imply that $v$ belongs to the set
\[
H \left(\varrho/2; \delta, \Lambda, \gamma \right)
\]
 of approximatively weakly symmetrical caloric functions. Consequently, using the symmetric caloric approximation lemma, i.e. Lemma \ref{lem5.1},  we the obtain existence of a caloric function $h$ that locally approximates $v$; more precisely
\begin{equation}\label{6.1.5}
 \exists \;  h\in H \left(\varrho/4;  \Lambda, \gamma \right) \quad \text{ such that} \quad  \dashint_{Q_{\varrho/4}(z_0)} \left| \frac{h-v}{\varrho/4} \right|^2  +  \gamma^{p-2}\left| \frac{h-v}{\varrho/4} \right|^p  \le \varepsilon
\end{equation}
Having such approximation of $v$ by $h$, we are ready  to show \eqref{shrink_step}; to this end, estimate $\psi_{s , z_0, l^{(2)}_{{\sigma\varrho}}} ({\sigma\varrho}) $, (which for $s$ being  $ 2, p $ constitute by definition $ \tilde  E_{z_0, l^{(2)}_{{\sigma\varrho}}} \left(\sigma\varrho \right)$) as follows
\begin{multline}\label{6.1.10}
\left({\theta\varrho\over4}\right)^{-s}\diagintop_{Q_{\varrho\theta\over4}(z_0)}\left|u-l^{(2)}_{{\theta\varrho\over4}}\right|^s\le 2^{s-1}\left({\theta\varrho\over4}\right)^{-s} \left[\diagintop_{Q_{\varrho\theta\over4}(z_0)}\left|u-l^{(s)}_{{\theta\varrho\over4}}\right|^s + \left|l^{(2)}_{{\theta\varrho\over4}}-l^{(s)}_{{\theta\varrho\over4}}\right|^s \right] \le 2^{s} C_\eqref{1.2.4} (n,s) \left({\theta\varrho\over4}\right)^{-s} \diagintop_{Q_{\varrho\theta\over4}(z_0)}\left|u-l^{(s)}_{{\theta\varrho\over4}}\right|^s  \\
\le
2^{s} C_\eqref{1.2.4} (n, s) \left({\theta\varrho\over4}\right)^{-s}  \diagintop_{Q_{\varrho\theta\over4}(z_0)}\left|u-l^{(2)}_\varrho -\gamma
[(h)_{z_0,{\theta\varrho\over4}}-(\nabla h)_{z_0,{\varrho\theta\over4}}(x-x_0)] \right|^s  \\
= 2^{s} C_\eqref{1.2.4} (n,s)   \gamma^s \left({\theta\varrho\over 4}\right)^{-s} \diagintop_{Q_{\varrho\theta\over4}(z_0)}|v-(h)_{z_0,{\theta\varrho\over4}}- (\nabla h)_{z_0,{\varrho\theta\over4}}(x-x_0)|^s \\
\le  2^{2s-1} C_\eqref{1.2.4} (n, s)   \gamma^2 \left[ \theta^{-(s+ n+2)}  \gamma^{s-2}  \dashint_{Q_{\varrho/4}(z_0)}\left| \frac{h-v}{\varrho/4} \right|^s+
  \left({\theta\varrho\over 4}\right)^{-s}   \gamma^{s-2}   \dashint_{Q_{\varrho\theta\over4}(z_0)}|h-(h)_{z_0,{\theta\varrho\over4}}- (\nabla h)_{z_0,{\varrho\theta\over4}}(x-x_0)|^s \right],
\end{multline}
where the second inequality holds in view of \eqref{1.2.4} of Lemma \ref{lem1.2}, the third one by minimization property of $l^{(s)}$ and the equality is given by definition \eqref{def_v} of $v$. To proceed further denote the mean integral over space (emphasizing its time dependance) by
\[
(g)(t)_{x_0,{\varrho}} \equiv \dashint_{B_{\varrho}(x_0)} g (x,t) dx
\] and estimate the second integral in the r.h.s. of \eqref{6.1.10} as follows
\begin{multline} \label{6.1.11}
 \dashint_{Q_{\varrho\theta\over4}(z_0)} \left|h-(h)_{z_0,{\theta\varrho\over4}}- (\nabla h)_{z_0,{\varrho\theta\over4}}(x-x_0) \right|^s \le \\
2^{s-1} \diagintop_{Q_{\varrho\theta\over4}(z_0)} \left|h(x,t)-(h)_{x_0,{\theta\varrho\over4}} (t)- (\nabla h)_{z_0,{\varrho\theta\over4}} (x-x_0) \right|^s dxdt+ 2^{s-1} \diagintop_{Q_{\varrho\theta\over4}(z_0)} \left| (h)_{x_0,{\varrho\theta\over4}} (t)- ( h)_{z_0,{\varrho\theta\over4}} \right|^s dxdt
\end{multline}
Observe that we cannot take $(h)_{x_0,{\varrho\theta\over4}} (t)$ instead of $( h)_{z_0,{\varrho\theta\over4}}$ directly in the second inequality of \eqref{6.1.10}, as only time-independent affine functions are admissible there.\\
 Consider the right-hand-side of  \eqref{6.1.11}.
For every $t$ one has  $(h-(h)_{x_0,{\theta\varrho\over4}} (t)- (\nabla h)_{z_0,{\varrho\theta\over4}}(x-x_0))_{x_0,{\theta\varrho\over4}}=0$, so Poincar\'e inequality in space followed by integration over time gives
\begin{multline}\label{6.1.115}
 \diagintop_{Q_{\varrho\theta\over4}(z_0)} \left|h-(h)_{x_0,{\theta\varrho\over4}} (t)- (\nabla h)_{z_0,{\varrho\theta\over4}} (x-x_0) \right|^s  \le \left({\theta\varrho\over 4}\right)^{s} \diagintop_{Q_{\varrho\theta\over4}(z_0)} \left| \nabla h - (\nabla h)_{z_0,{\varrho\theta\over4}} \right|^s \le \\
  C_\eqref{2.4.12} (\lambda, \Lambda, s)  \theta^s  \left({\theta\varrho\over 4}\right)^{s} \diagintop_{Q_{\varrho\over4}(z_0)} | \nabla h|^s \le  C_\eqref{2.4.12} (\lambda, \Lambda, s)  \theta^s  \left({\theta \varrho\over 4}\right)^{s} \gamma^{2-s};
\end{multline}
the second inequality results from estimate \eqref{2.4.12} of Lemma \ref{LinSys} used with $m=0, \;q=p=s$ for $\frac{\partial h}{\partial_{x_i}}$;  the last inequality is valid as $h$ is a symmetrical caloric function. Simultaneously we have thanks to a smoothness of $h$ and the mean-value property
\begin{multline} \label{6.1.12}
\diagintop_{Q_{\varrho\theta\over4}(z_0)} \left| (h)_{x_0,{\varrho\theta\over4}} (t)- ( h)_{z_0,{\varrho\theta\over4}} \right|^s = \diagintop^{t_0}_{t_0 - \left({\varrho\theta\over4} \right)^2}   \left|  \diagintop^{t_0}_{t_0 - \left({\varrho\theta\over4} \right)^2}   (h)_{x_0,{\varrho\theta\over4}} (t)- ( h)_{x_0,{\varrho\theta\over4}} (r) dr\right|^s   dt  \le  \\
\diagintop^{t_0}_{t_0 - \left({\varrho\theta\over4} \right)^2} \left|\;   \diagintop^{t_0}_{t_0 - \left({\varrho\theta\over4} \right)^2} \left| \sup_{ \tau \in \left({t_0 - \left({\varrho\theta\over4} \right)^2}, {t_0} \right)} \left(\frac{\partial h}{\partial \tau} \right)_{x_0,{\varrho\theta\over4}} (\tau) \right| |t-r|d r \right|^s   dt  \le  \left({\varrho\theta\over4} \right)^{2s} \left| \sup_{ \tau \in \left({t_0 - \left({\varrho\theta\over4} \right)^2}, {t_0} \right)} \left(\frac{\partial h}{\partial \tau} \right)_{x_0,{\varrho\theta\over4}} (\tau) \right|
 \end{multline}
Using the inequality \eqref{2.4.1} of Lemma \ref{LinSys} with $m=1, \;q= \infty, \;p=s$ to estimate the r.h.s. of \eqref{6.1.12} one arrives at
\begin{equation} \label{6.1.13}
\diagintop_{Q_{\varrho\theta\over4}(z_0)} \left| (h)_{x_0,{\varrho\theta\over4}} (t)- ( h)_{z_0,{\varrho\theta\over4}} \right|^s \le  C_\eqref{2.4.1} (\lambda,  \Lambda, s)   \left({\varrho\theta\over4} \right)^{2s} \left(\varrho\over4 \right)^{-s}\diagintop_{Q_{\varrho\over4}(z_0)} \left|\frac{h}{\varrho / 4} \right|^s \le C_\eqref{2.4.1} (\lambda,  \Lambda, s)  \theta^s \left({\varrho\theta\over4} \right)^{s} \gamma^{2-s}
\end{equation}
where the second inequality results from $h$ being a symmetrical caloric function. Combine \eqref{6.1.115} and \eqref{6.1.13} to estimate the right-hand-side of \eqref{6.1.11}
\begin{equation} \label{6.1.11'}
 \dashint_{Q_{\varrho\theta\over4}(z_0)} \left|h-(h)_{z_0,{\theta\varrho\over4}}- (\nabla h)_{z_0,{\varrho\theta\over4}}(x-x_0) \right|^s \le
2^{s-1} \left( C_\eqref{2.4.1} (\lambda,  \Lambda, s) + C_\eqref{2.4.12} (\lambda,  \Lambda, s)\right)  \theta^s \left({\varrho\theta\over4} \right)^{s} \gamma^{2-s}
\end{equation}
This and \eqref{6.1.5} we use in \eqref{6.1.10} to get
\begin{equation}\label{6.1.10'}
\left({\theta\varrho\over4}\right)^{-s}\diagintop_{Q_{\varrho\theta\over4}(z_0)}\left|u-l^{(2)}_{{\theta\varrho\over4}}\right|^s 
\le  2^{2s-1}  C_\eqref{1.2.4} (n,s) \gamma^2 \left[ \theta^{-(s+ n+2)}  \varepsilon+ \theta^s
2^{s-1} \left( C_\eqref{2.4.1} (\lambda,  \Lambda, s) + C_\eqref{2.4.12} (\lambda,  \Lambda, s)\right) \right]
\end{equation}
Estimate \eqref{6.1.10'} with $s=2$ and $s=p$ gives, in view of $\theta \le 1$, Definition \ref{energies} of $\tilde E$ and choice \eqref {choice_gamma} of $\gamma$
\begin{equation}\label{6.1.56'} 
E_{z_0, l^{(2)}_{{\theta\varrho\over4}}} \left( {\theta\varrho\over4} \right) \le  C^2_{\ref{cor4.2}}(M)
 \left[E_{z_0, l^{(2)}_\varrho} (\varrho) + (\delta_\varepsilon/2)^{-2}\varrho^{2\beta}\right]  C_\eqref{6.1.56'}   (\lambda,  \Lambda, p) \left[ \theta^{-(p+ n+2)}  \varepsilon+ \theta^2 \right] 
\end{equation}
where (robustly)
\begin{equation}
C_\eqref{6.1.56'}   \equiv  2^{5p -3}\max_{s \in \{2;p\}}  C_\eqref{1.2.4} (n,s)   \left( C_\eqref{2.4.1} (\lambda,  \Lambda, s) + C_\eqref{2.4.12} (\lambda,  \Lambda, s)\right) 
\end{equation}
recall that we have taken $\sigma = \theta/4$; this with the definition of $C_\eqref{6.1.05const} $ gives from \eqref{6.1.56'} 
\begin{equation}\label{6.1.56''} 
 \tilde  E_{z_0, l^{(2)}_{{\sigma \varrho}}} \left( {\sigma \varrho} \right) \le  C_\eqref{6.1.05const}  \left[ \frac{ \varepsilon}{ (4 \sigma)^{p+ n+2}}+ 16 \sigma^2 \right]  \left[E_{z_0, l^{(2)}_\varrho} (\varrho) + (\delta_\varepsilon/2)^{-2}\varrho^{2\beta}\right]  +  ({\sigma \varrho})^{2 \beta} \le C_\eqref{6.1.05const} 32 \sigma^2 4 \left[E_{z_0, l^{(2)}_\varrho} (\varrho) + \delta_\varepsilon^{-2}\varrho^{2\beta}\right]  +  ({\sigma \varrho})^{2 \beta}
\end{equation}
where the second inequality is given by choice of $\varepsilon $, see \eqref{6.1.05}; in the same inequality we have chosen $\sigma$ so that $C_\eqref{6.1.05const} 2^7 \sigma^2 \le  \sigma^{2 \alpha}$, which gives the thesis.
\end{proof}
Let us now state inequalities used for the singular set description in the following iteration of Lemma \ref{lem6.1} performed in Lemma \ref{lem6.2}.
\begin{lem}\label{do_iter}
Under the assumptions of Lemma \ref{lem6.1}, the following inequalities hold
\begin{multline}\label{pom_itern}
\varrho^{-p} \diagintop_{Q_\varrho(z_0)}|u-l^{(2)}_{\varrho} |^p \le \\
C(|(u)_{z_0}| + |(\nabla u)_{z_0}|)  \left[  \varrho^{p\beta} + \max_{l=1,\beta} \left(  
 \diagintop_{Q_\varrho(z_0)} | \nabla u - (\nabla u)_{z_0}|^p  \right)^l+
   \max_{r=1,\beta(p-1),p-1}   \left( \diagintop_{Q_\varrho(z_0)} |\D u-\ (\D u)_{z_0}  |^p  \right)^{r}  \right]\end{multline}
\begin{multline}\label{pom_iter}
\varrho^{-p} \diagintop_{Q_\varrho(z_0)}|u-(u)_{z_0} - (\nabla u)_{z_0} (x-x_0) |^p \le \\
C(|(u)_{z_0}| + |(\nabla u)_{z_0}|)  \left[  \varrho^{p\beta} + \max_{l=1,\beta} \left(  
 \diagintop_{Q_\varrho(z_0)} | \nabla u - (\nabla u)_{z_0}|^p  \right)^l+
   \max_{r=1,\beta(p-1),p-1}   \left( \diagintop_{Q_\varrho(z_0)} |\D u-\ (\D u)_{z_0}  |^p  \right)^{r}  \right]\end{multline}
\begin{multline}\label{pom_iter'}
\varrho^{-p}  \diagintop_{Q_\varrho(z_0)} |u- (u)_{z_0} |^p \le C(|(u)_{z_0}| + |(\nabla u)_{z_0}|)  \left( 1+   \left( \diagintop_{Q_\varrho(z_0)} |\D u-\ (\D u)_{z_0}  |^p  \right)^{p-1} + \diagintop_{Q_\varrho(z_0)} | \nabla u - (\nabla u)_{z_0}|^p   \right)
\end{multline}
\end{lem}
\begin{proof}
Recall that  $\eta_\varrho (x)$ denotes a standard space-mollifier, supported in $B_\varrho(x_0)$. Define
\begin{equation}
(u)_{\eta, x_0} (t) \equiv \int_{B_\varrho(z_0)} u(x,t) \eta_\varrho (x) dx, \qquad (u)_{\eta, z_0} \equiv \dashint_{t_0 - \varrho^2}^{t_0} (u)_{\eta, x_0} (t)
\end{equation}
Let us perform estimates, using first \eqref{1.2.4}, then the approximative minimization property of standard mean value with respect to $L^p$ norms
\begin{multline}\label{6.3.1to'}
 \diagintop_{Q_\varrho(z_0)}|u-(u)_{z_0} - (\nabla u)_{z_0} (x-x_0) |^p \le  C \diagintop_{Q_\varrho(z_0)} \left|u- (u)_{\eta, z_0} - (\nabla u)_{z_0} (x-x_0) \right|^p \le \\
C \left [\diagintop_{Q_\varrho(z_0)} \left|u- (u)_{\eta, x_0} (t) - (\nabla u)_{z_0} (x-x_0) \right|^p dxdt+ \diagintop_{Q_\varrho(z_0)} | (u)_{\eta, x_0} (t) - (u)_{\eta, z_0} |^p dxdt \right]
\end{multline}
The first integral in the right-hand-side of \eqref{6.3.1to'} is majorized in view of Poincar\'e inequality in space (for generalized integral means) by
\begin{equation}\label{6.3.1toto1'}
C \varrho^{p} \diagintop_{Q_\varrho(z_0)} |\nabla u - (\nabla u)_{z_0} |^p 
\end{equation}
whereas for the second one, in view of the inequality \eqref{pre4.1.1} of Lemma \ref{lem4.1} holds
\begin{multline}\label{nwlem_1}
\diagintop_{Q_\varrho(z_0)} | (u)_{\eta, x_0} (t) - (u)_{\eta, z_0}|^p = \diagintop^{t_0}_{t_0 -\varrho^2} \left|  \diagintop^{t_0}_{t_0 -\varrho^2}  (u)_{\eta, x_0} (t) - (u)_{\eta, x_0} (r) dr \right|^p dt \le \diagintop^{t_0}_{t_0 -\varrho^2} \diagintop^{t_0}_{t_0 -\varrho^2}  \left| \intop_{B_\varrho(x_0)} (u (t,x) - u(r,x) )\eta_\varrho(x) \right|^p drdt \\
\le C^p_\eqref{pre4.1.1}  ( |(u)_{z_0}| + |(\D u)_{z_0}|) \varrho^p \left|  \diagintop_{Q_\varrho(z_0)}  \varrho^\beta + |u- (u)_{z_0} |^\beta+ |u- (u)_{z_0} | + |\D u-\ (\D u)_{z_0} |^{p-1} + |\D u-\ (\D u)_{z_0} | \right|^p \le \\
C \cdot C^p_\eqref{pre4.1.1}  (|(u)_{z_0}| + |(\D u)_{z_0}|) \varrho^p  \left[  \varrho^{p\beta} + \left( \diagintop_{Q_\varrho(z_0)} |u- (u)_{z_0} |^p  \right)^\beta+ \diagintop_{Q_\varrho(z_0)} |u- (u)_{z_0} |^p  + \max_{r=1,p-1}   \left( \diagintop_{Q_\varrho(z_0)} |\D u-\ (\D u)_{z_0}  |^p  \right)^{r}  \right]
\end{multline}
Summing up, \eqref{6.3.1to'} takes the form
\begin{multline}\label{nwlem_3}
\diagintop_{Q_\varrho(z_0)}|u-(u)_{z_0} - (\nabla u)_{z_0} (x-x_0) |^p \le
C \varrho^p\diagintop_{Q_\varrho(z_0)} |\nabla u - (\nabla u)_{z_0} |^p +\\
C \varrho^p C^p_\eqref{pre4.1.1}  (|(u)_{z_0}| + |(\D u)_{z_0}|) \left[  \varrho^{p\beta} + \max_{l=1,\beta} \left( \diagintop_{Q_\varrho(z_0)} |u- (u)_{z_0} |^p  \right)^l+   \max_{r=1,p-1}   \left( \diagintop_{Q_\varrho(z_0)} |\D u-\ (\D u)_{z_0}  |^p  \right)^{r}  \right]\end{multline}

Further, estimate 
\begin{multline}\label{nwlem_2}
\diagintop_{Q_\varrho(z_0)} |u- (u)_{z_0} |^p \le C \diagintop_{Q_\varrho(z_0)} \left[ |u- (u)_{\eta, x_0} (t) |^p + |(u)_{\eta, x_0} (t)-(u)_{\eta, z_0} |^p\right]  dxdt + C \diagintop^{t_0}_{t_0 -\varrho^2}   | (u)_{\eta, x_0} (t) - (u)_{x_0} (t)|^p dt \le \\
C\varrho^p \diagintop_{Q_\varrho(z_0)} | \nabla u|^p + C \diagintop^{t_0}_{t_0 -\varrho^2} \diagintop^{t_0}_{t_0 -\varrho^2}  \left| \intop_{B_\varrho(x_0)} (u (t,x) - u(r,x) )\eta_\varrho(x) dx \right|^p drdt  + \\
 C \diagintop_{Q_\varrho(z_0)}  | u - (u)_{\eta, x_0} (t) - (\nabla u)_{z_0} (x-x_0) |^p + | u - (u)_{x_0} (t) - (\nabla u)_{z_0} (x-x_0)|^p dt \\
\le C\varrho^p \left( \diagintop_{Q_\varrho(z_0)} | \nabla u|^p +C^p_{\eqref{prepre4.1.1} } (|(u)_{z_0}| + |(\D u)_{z_0}|)  \left[\diagintop_{Q_\varrho(z_0)}  1 + |\D u|^{p-1} \right]^p + \diagintop_{Q_\varrho(z_0)} | \nabla u - (\nabla u)_{z_0}|^p \right) \\
 \le
C(|(u)_{z_0}| + |(\nabla u)_{z_0}|)  \varrho^p \left(1+  \left( \diagintop_{Q_\varrho(z_0)} |\D u-\ (\D u)_{z_0}  |^p  \right)^{p-1} + \diagintop_{Q_\varrho(z_0)} | \nabla u - (\nabla u)_{z_0}|^p   \right)
\end{multline}
where for the second inequality we use Poincar\'e in space and computation analogous to that of \eqref{nwlem_1} and for the third one: estimate \eqref{prepre4.1.1} and again Poincar\'e (both for standard and generalized integral means). This inequality is  \eqref{pom_iter'}; it used in \eqref{nwlem_3} gives
\begin{multline}\label{nwlem_3'}
\varrho^{-p}\diagintop_{Q_\varrho(z_0)}|u-(u)_{z_0} - (\nabla u)_{z_0} (x-x_0) |^p \le
C \diagintop_{Q_\varrho(z_0)} |\nabla u - (\nabla u)_{z_0} |^p +\\
C(|(u)_{z_0}| + |(\nabla u)_{z_0}|)  \left[  \varrho^{p\beta} + \max_{l=1,\beta} \left( \varrho^p+  
 \diagintop_{Q_\varrho(z_0)} | \nabla u - (\nabla u)_{z_0}|^p  \right)^l+
   \max_{r=1,\beta(p-1),p-1}   \left( \diagintop_{Q_\varrho(z_0)} |\D u-\ (\D u)_{z_0}  |^p  \right)^{r}  \right]\end{multline}
which gives \eqref{pom_iter}.  Finally to obtain \eqref{pom_itern} from \eqref{pom_iter}, estimate from below l.h.s. of \eqref{6.3.1to'}, using first \eqref{1.2.4}, then the minimization property of $l^{(p)}_{\varrho} $
\begin{equation}\label{6.3.1to''}
\diagintop_{Q_\varrho(z_0)} \left|u-l^{(2)}_{\varrho} \right|^p \le 2^{p-1}(C_\eqref{1.2.4}(n,p) + 1) \diagintop_{Q_\varrho(z_0)} \left|u-l^{(p)}_{\varrho} \right|^p \le C \diagintop_{Q_\varrho(z_0)} \left|u- (u)_{\eta, z_0} - (\nabla u)_{z_0} (x-x_0) \right|^p 
\end{equation}

\end{proof}

\begin{lem}\label{lem6.2} Let $p \ge 2$ and $u \in C(-T, 0; L^2(\Omega)) \cap L^p(-T,0; W^{1,p}(\Omega))$  be a weak solution to (\ref{1.1}) under structure conditions (\ref{1.2a} --- \ref{1.5}). Take $z_0\in Q_T$ such that
\begin{equation}\label{6.3.1}
\liminf_{\varrho\to0}   \diagintop_{Q_\varrho(z_0)} |\nabla u- (\nabla u)_{z_0}  |^p  =0,
\end{equation}
\begin{equation}\label{6.3.2}
\limsup_{\varrho\to0}|(u)_{z_0,\varrho}|+|(\nabla u)_{z_0,\varrho}|<+\infty
\end{equation}
then in $\tilde Q (z_0)$, denoting a certain vicinity of $z_0$, holds
\begin{equation}
\nabla u\in C^{\beta,{\beta\over 2}}(\tilde Q (z_0))
\end{equation}
where $ \beta$ is given by (\ref{1.8}).
\end{lem}
\begin{proof}
As in view of Lemma \ref{lem1.2} $( u)_{z_0, \varrho} =  l^{(2)}_{\varrho} (x_0) $, assumptions \eqref{6.3.1}, \eqref{6.3.2} and pointwise estimate $|\D g| \le |\nabla g| $ imply that we can find sequence $\varrho_n \to 0$ for which the following hold for a certain $M < \infty$ 
\begin{equation}\label{p6.3.45}
\lim_{{\varrho_n}\to0}  \left[ \diagintop_{Q_{\varrho_n}(z_0)} |\D u-\ (\D u)_{z_0}  |^p  +   \diagintop_{Q_{\varrho_n}(z_0)} |\nabla u- (\nabla u)_{z_0}  |^p \right]  =0
\end{equation}
\begin{equation}\label{6.3.45}
|(\nabla u)_{z_0,\varrho_n} |\le {M/ 8}, \quad | l^{(2)}_{\varrho_n} (x_0) |\le{M/ 8}, \quad  \tilde  E_{l^{(2)}_{{\varrho_n}}}  \to 0
\end{equation}
where H\"older inequality and $p \ge 2$ is used to control $\psi_2$ in $ \tilde  E $ with $\psi_p$ being l.h.s. of  \eqref{pom_itern}. Moreover  \eqref{p6.3.45}, \eqref{6.3.45} with inequality \eqref{pom_iter} give again via H\"older inequality 
\begin{equation}\label{6.3.455}
\lim_{{\varrho_n}\to0}  \varrho_n^{-2}   \diagintop_{Q_{\varrho_n}(z_0)}| u-(u)_{z_0,{\varrho_n}} - (\nabla u)_{z_0,{\varrho_n}} (x-x_0) |^2 = 0
\end{equation}
 In order to replace $|(\nabla u)_{z_0,\varrho} |$  in \eqref{6.3.45} with $\left|\nabla l^{(2)}_{z_0, \varrho_n}  \right|$, perform estimate using inequality \eqref{1.2.2} of Lemma \ref{lem1.2}
\begin{equation}\label{6.3.47}
|\nabla l^{(2)}_{z_0, \varrho_n}  |^2 \le {2d(d+2)\over{\varrho_n}^2}   \diagintop_{Q_{\varrho_n}(z_0)}| u-(u)_{z_0,{\varrho_n}} - (\nabla u)_{z_0,{\varrho_n}} (x-x_0) |^2 + |(\nabla u)_{z_0,\varrho_n}|^2,
\end{equation} 
which in view of  \eqref{6.3.455} and  \eqref{6.3.47} gives for $n \ge n_0$
\begin{equation}
|\nabla l^{(2)}_{z_0, \varrho_n}  | \le M/4
\end{equation}
This and \eqref{6.3.45} shows that there is a sequence $\varrho_m \to 0$ ($m = n- n_0$) for which holds for a certain $M < \infty$
\begin{equation}\label{6.3.5}
 |\nabla l^{(2)}_{z_0, \varrho_m}  |\le{M/4}, \quad | l^{(2)}_{\varrho_m}  (x_0) |\le{M/ 4}, \quad  \tilde  E_{l^{(2)}_{{\varrho_m}}}  \to 0
\end{equation}
Constant $M$ and choice  $\alpha \in ( \beta, 1) $, fixes $\sigma \in (0, 1/4), \; \delta \in (0,1)$ as in Lemma \ref{lem6.1}. Convergences in \eqref{6.3.5} imply the existence of such $\varrho_0$ that
\begin{equation}\label{p6.3.53}
\omega(M+1, \tilde  E_{ l^{(2)}_{z_0, \varrho_0}} ({\varrho_0}) )+ \tilde  E^{\frac{1}{2}}_{l^{(2)}_{z_0, \varrho_0}} ({\varrho_0})  <{\delta \over 2 C_{\ref{cor4.2}}(M) }
\end{equation}
holds; in fact, for this fixed radius $\varrho_0$, the absolute continuity of integrals with respect to the Lebesgue measure and a continuity of the modulus of continuity $\omega$ imply that we have for any point of $\tilde z \in \tilde Q (z_0)$, $\tilde Q (z_0)$ being a certain neighborhood of $z_0$,
\begin{equation}\label{pp6.3.53}
 |\nabla l^{(2)}_{{\tilde z}, \varrho_m}  |\le{M/4}, \quad | l^{(2)}_{{\tilde z}, \varrho_m}  ({\tilde x}) |\le{M/4}, \quad 
\omega(M+1, \tilde  E_{ l^{(2)}_{{\tilde z}, \varrho_0}} ({\varrho_0}) )+ \tilde  E^{\frac{1}{2}}_{l^{(2)}_{{\tilde z}, \varrho_0}} ({\varrho_0})  <{\delta \over 2 C_{\ref{cor4.2}}(M) }
\end{equation}
in what follows, we generally abandon the dependence of the following expressions on $\tilde z$, remembering that it is an arbitrary point from $ \tilde Q (z_0)$. Lemma \ref{lem6.1} and  \eqref{pp6.3.53} give
\begin{equation}\label{6.3.54}
 \tilde  E_{ l^{(2)}_{{\sigma \varrho_0}}} \left(\sigma \varrho_0 \right) \le  \sigma^{2\beta}  \tilde  E_{l^{(2)}_{\varrho_0}} (\varrho_0) + C(M, \beta) (\sigma \varrho_0)^{2\beta}
 \end{equation}
The next step is to prove that for every $j \in \en$ holds
\begin{equation}\label{6.3.6}
 \tilde  E_{ l^{(2)}_{{\sigma^j \varrho_0}}} \left(\sigma^j \varrho_0 \right) \le  \sigma^{j 2\beta}  \tilde  E_{l^{(2)}_{\varrho_0}} (\varrho_0) +  C (M, \beta)  (\sigma^j \varrho_0)^{2\beta}, \qquad 
|l^{(2)}_{\sigma^{j}\varrho_0}|+|\nabla l^{(2)}_{\sigma^{j}\varrho_0}|\le M
\end{equation}
as the inductive argument here is identical as the respective part of the proof of Lemma 4.9 in \cite{[DMS]}, we do not present it here.\\
Now we show, for the above fixed $\varrho_0$, $\sigma$, $M$, that
\begin{equation}
\label{6.3.67D}
\lim_{j\to\infty}(\nabla u)_{\tilde z,{\sigma^j\varrho_0}} \equiv \tilde \Gamma
\end{equation}
exists and for $r\in\left(0,{\varrho_0 /2}\right)$
\begin{equation}
\label{6.3.7D}
\diagintop_{Q_r(\tilde z)}|\nabla u-\tilde  \Gamma|^2\le C  r^{2\beta}
\end{equation}
Fix $r\in\left(0,{\varrho_0 /2}\right)$ and choose $j$ such that
\begin{equation}\label{6.3.75}
\sigma^{j+1}{(\varrho_0/2)}<r\le\sigma^{j}{(\varrho_0/2)}
\end{equation}
Then by the minimizing property of a mean value (the first inequality), the Caccioppoli inequality \eqref{4.1.1'} (the middle inequality),  \eqref{6.3.6} (the third one) we get
\begin{multline}\label{6.3.9}
 \diagintop_{Q_r} |\nabla u-(\nabla u)_{Q_r}|^2 \le 
  (1 / \sigma )^{n+2} \diagintop_{Q_{\sigma^{j}(\varrho_0/2)}} \left|\nabla u-\nabla l^{(2)}_{\sigma^{j}(\varrho_0/2)} \right|^2  \le 
  (1 / \sigma )^{n+2}    \tilde  E_{ l^{(2)}_{{\sigma^j \varrho_0}}} \left(\sigma^j \varrho_0 \right) \le \\
   (1 / \sigma )^{n+2}   \left[   \sigma^{j 2\beta}  \tilde  E_{l^{(2)}_{\varrho_0}} (\varrho_0) +  C(M, \beta)  (\sigma^j \varrho_0)^{2\beta} \right] \le     (1 / \sigma )^{n+2} \left[  \tilde  E_{l^{(2)}_{\varrho_0}} (\varrho_0) +  C (M, \beta)   \varrho_0^{2\beta} \right] \sigma^{j 2\beta} \le C({M, \beta})r^{2\beta}
\end{multline}
Similarly for $j< k$
\begin{equation}\label{6.3.10}
\begin{aligned}
&\left|(\nabla u)_{\sigma^j(\varrho/2)}-(\nabla u)_{\sigma^k(\varrho/2)} \right|\le\sum_{m=j+1}^k
\left| (\nabla u)_{\sigma^m(\varrho/2)}-(\nabla u)_{\sigma^{m-1}(\varrho/2)} \right|\\
&\le\sum_{m=j+1}^k\left|\diagintop_{Q_{\sigma^m(\varrho/2)}}\nabla u-(\nabla u)_{\sigma^{m-1}(\varrho/2)}\right|\le
  (1 / \sigma )^{\frac{n+2}{2}}  \sum_{m=j+1}^k \left[\diagintop_{Q_{\sigma^{m-1}(\varrho/2)}}
\left|\nabla u-(\nabla u)_{\sigma^{m-1}(\varrho/2)} \right|^2\right]^{1\over2}\\
&\le   (1 / \sigma )^{\frac{n+2}{2}}    \left[  \tilde  E_{l^{(2)}_{\varrho_0}} (\varrho_0) +  C(M, \beta)   \varrho_0^{2\beta} \right]^{1/2}  \sum_{m=j}^{k-1}\sigma^{m \beta} \le    (1 / \sigma )^{\frac{n+2}{2}}  (1- \sigma)^{-1} C^{1/2} ({M, \beta}) \sigma^{\beta j} = \tilde C ({M, \beta})\sigma^{\beta j}
\end{aligned}
\end{equation}
where the last but one inequality is obtained as (\ref{6.3.9}). The estimate (\ref{6.3.10}) states that 
$\lim_{j\to\infty}(\nabla u)_{\sigma^j{\varrho\over2}}=\tilde \Gamma$ exists and that
\begin{equation}\label{6.3.11}
|(\nabla u)_{\tilde z, \sigma^j(\varrho/2)}- \tilde \Gamma|\le  \tilde C ({M, \beta}) \sigma^{\beta j}
\end{equation}
This combined with (\ref{6.3.9}) results in:
\begin{equation}\label{6.3.12}
\diagintop_{Q_r}|\nabla u- \tilde \Gamma|^2\le    (1 / \sigma )^{n+2} \diagintop_{Q_{\sigma^{j}(\varrho_0/2)}{} }
|\nabla u\pm(\nabla u)_{\sigma^{j}(\varrho_0/2)} -\tilde \Gamma|^2\le2 (  \tilde C^2 ({M, \beta}) + C ({M, \beta})) \sigma^{2\beta j} \le  C ({M, \beta}) r^{2\beta}
\end{equation}
where the last inequality holds in view of  \eqref{6.3.75}.
As (\ref{6.3.12}) is valid for any $\tilde z$, being an arbitrary point from $ \tilde Q (z_0)$, imbedding of Campanato into H\"older spaces gives
\begin{equation*}
\nabla u\in C^{\beta,{\beta\over 2}}(\tilde Q (z_0))
\end{equation*}
\end{proof}
We are done with the partial regularity result for the gradient. Let us now focus on an analogous property for the solution itself, stated in the following result.
\begin{lem}\label{lem6.2'} Take $ \tilde Q ( z_0 )$ -- a neighborhood of a regular point of Lemma \ref{lem6.2}. Under assumptions of Lemma \ref{lem6.2} holds
\begin{equation}
 u\in C^{1,{1\over 2}}(  \tilde Q ( z_0 )).
\end{equation}
\end{lem}
\begin{proof}
For any $ x,y \in B_{\tilde \varrho} (\tilde x) \subset \tilde Q ( z_0 )$ we have a pointwise estimate for a 
$ C(-T, 0; L^2(\Omega)) \cap L^p(-T,0; W^{1,p}(\Omega))$ weak solution to (\ref{1.1})

\begin{equation}\label{bojhaj}
\left| \intop_{ B_{\tilde \varrho} (\tilde x)} (u (\tilde \tau ,x) - u(\tilde \tau,y) )\eta_\varrho(x) dx \right| \le C \diagintop_{ B_{\tilde \varrho} (\tilde x)} (M (\nabla u) (t,x)+ M (\nabla u) (t,y) ) |x-y| dx  \le C_{ \tilde Q ( z_0 )} \tilde \varrho,
\end{equation}
where $M(f)$ is a maximal function.  The first inequality is given by Bojarski-Haj\l asz inequality (see \cite{[BH]}, Theorem 3) and the second is a consequence of boundedness of gradients given by  Lemma \ref{lem6.2}. Adding \eqref{prepre4.1.1}, which holds for every time level, and twice \eqref{bojhaj}, one with $\tilde \tau = \tau, y = y_1$ and the second with $\tilde \tau =t,  y = y_2$ we obtain thanks to boundedness of $\nabla u$
\begin{multline}\label{bd_u}
C_{ \tilde Q ( z_0 )} \tilde \varrho \ge \\
\left| \intop_{B_{\tilde \varrho} (\tilde x)} (u (t,x) - u(\tau,x) )\eta_\varrho(x) dx \right|  + \left| \intop_{ B_{\tilde \varrho} (\tilde x)} (u (t ,x) - u(t ,y_1) )\eta_\varrho(x) dx \right| + \left| \intop_{ B_{\tilde \varrho} (\tilde x)} (u (\tau ,x) - u(\tau ,y_2) )\eta_\varrho(x) dx \right|  \ge |   u(t ,y_1) - u(\tau ,y_2)| \end{multline}
which gives thesis.
\end{proof}
Finally we see that
\begin{proof}[Proof of Theorem \ref{theo1.1}]
results from Lemmas \ref{lem6.2} and \ref{lem6.2'}.
\end{proof}
\section{Conclusions}
\setcounter{equation}{0}
The natural next step is to perform an analysis of the Hausdorff dimension of the singular set, at least for less general systems, for example for which the dependence of the main part $A$ on $u$ is waived. This, together with the non-linear Calderon-Zygmund $L^q$ estimates will be the joint content of the forthcoming paper, as there is a natural connection between the singular set estimates and the restriction on $q$.\\
 It would be interesting, using new results on parabolic approximation, to perform similar analysis for $p$-Stokes system. Finally let us mention, that it seems that for a certain range of $p$'s, close to $2$, full $C^{1, \alpha}$ regularity for symmetric $p$-Laplace holds; this is also currently work in progress. \\However, the ultimate goal in this field, namely the full interior $C^{1, \alpha}$-regularity for symmetric $p$-Laplace system, without restrictions on $p$ and the space dimensions, seems to be essentially open.

\section{Appendix} \label{app}
\setcounter{equation}{0}
Here we present results which have been removed from the main part of this article for the sake of  traceability.
\begin{proof}[Proof of Korn's inequality \eqref{1.3.2} in Lemma \ref{lem1.3}]
Use inequality from \cite{[DV]}
\begin{equation}\label{1.3.3}
K\intop_{B_1(x)}|\D h|^2\ge\inf_{R\in{\cal R}}\intop_{B_1(x)}|\nabla[h-R]|^2
\end{equation}
where ${\cal R}$ is the set of rigid motions, i.e. affine functions with antisymmetric linear part. (\ref{1.3.3}) with $h:=u-(\D u)(x-x_0)$ yields:
\begin{equation*}
K\intop_{B_1(x)}|\D u-(\D u)|^2\ge\inf_{R\in{\cal R}}\intop_{B_1(x)}|\nabla u-(\D u)-R|^2\ge\intop_{B_1(x)}|\nabla u-(\nabla u)|^2
\end{equation*}
we have also pointwisely $|\nabla h|^2\ge|\D h|^2$, so for $h:=u-(\nabla u)(x-x_0)$ it gives
\begin{equation*}
\intop_{B_1(x)}|\nabla u-(\nabla u)|^2\ge\intop_{B_1(x)}|\D u-(\D u)|^2
\end{equation*}
The independence of $K$ on radius $r$ comes from scaling.
\end{proof}
Next we show the needed result on linear systems.
\begin{proof}[Proof of Lemma \ref{LinSys}]
 Smoothness is a standard result for systems with coefficients depending on full gradient and satisfying Legendre-Hadamard conditions. See \cite{[EZ]}, \cite{[Sch]}.
To prove inequalities we modify slightly the technique of Campanato \cite{[C]}. 
Scaling $v(y,s)=u(y/r,s/{r^2})$ justifies that $u$ solves locally \eqref{linC} in $Q_r$ iff $v$ solves  \eqref{linC}  locally in $Q_1$. Therefore we consider first $v$ satisfying
\begin{equation}\label{linCQ1}
\intop_{Q_1}v\varphi_{,t}- A \; \D v \; \D\phi =0 \quad \hforall_{\varphi\in C_0^\infty(Q_1)}
\end{equation}
  Take a smooth cutoff functions from  $C_0^\infty(Q_1)$ that satisfies
\begin{equation*}
[0,1]\ni\theta(x)=\begin{cases}1&B_{1\over2}\\ 0&B_1^c\end{cases} \quad |\nabla\theta|\le{8}
\end{equation*}
\begin{equation*}
[0,1]\ni \sigma(t)=\begin{cases}1&t\in \big({1\over4},0\big),\\ 
0&t\le 1\end{cases} \quad |\sigma_{,t}|\le 8
\end{equation*}
we test (\ref{2.4.1}) with $\theta^2\sigma^2v$, which yields:

\noindent
(i) for the main part: 
\begin{equation}\label{2.4.2}
\begin{aligned}
&{1\over2}\sigma^2a_{kl}^{ij}(v_l^k+v_k^l)(\theta^2v^i)_{,j}={1\over2}a_l^{ij}(v_l^k+v_k^l)[\theta(\theta v^i)_{,j}+\theta_{,j}\theta v^i]\\
&={\sigma^2\over 2} a_{kl}^{ij}[(v\theta)_{,l}^k+(v\theta)_{,k}^l](\theta v^i)_{,j}-{\sigma^2\over 2}a_{kl}^{ij}(v\theta_{,l}+v\theta_{,k})(\theta v^i)_{,j}+{\sigma^2\over 2} a_{kl}^{ij}(v_l^k+v_k^l)\theta v^i\theta_{,j}\\
&\ge\lambda|\D(v\sigma\theta)|^2+{\sigma^2\theta\over 2} [a_{kl}^{ij}(v_l^k+v_k^l)v^i\theta_{,j}-a_{kl}^{ij}(v^k\theta_{,l}+v^l\theta_{,k})v_{,j}^i]-|A||\sigma^2|\,|\nabla\theta|^2|v|^2
\end{aligned}
\end{equation}
and the middle part vanishes because $a_{kl}^{ij}=a_{ij}^{kl}=a_{lk}^{ji}$
\begin{equation*}
a_{kl}^{ij}v_l^kv^i\theta_{,j}+a_{kl}^{ij}v_k^lv^i\theta_{,j}-
\overbrace{a_{kl}^{ij}v_j^iv^k\theta_{,l}}^{=a_{ij}^{kl} v_l^kv^i\theta_{,j}}
-\overbrace{a_{kl}^{ij}v_j^iv^l\theta_{,k}}^{=a_{ji}^{lk}v_k^lv^i\theta_{,j}}
=\underbrace{(a_{kl}^{ij}-a_{ij}^{kl})}_{=0}v_{,l}^kv^i\theta_{,j}+\underbrace{(a_{kl}^{ij}-a_{lk}^{ji})}_{=0}
v_k^lv^i\theta_{,j}
\end{equation*}
(ii) for the evolutionary part:
\begin{equation}\label{2.4.3}
v_{,t}\theta^2\sigma^2v={1\over2}(v^2\theta^2\sigma^2)_{,t}-v^2\theta^2\sigma_{,t}\sigma
\end{equation}
By (\ref{2.4.2}) and (\ref{2.4.3}) integrated over $Q_1$ we have:
\begin{equation}\label{2.4.45}
{1\over2}{d\over dt}\intop_{Q_1}v^2\theta^2 \sigma^2+\lambda\intop_{Q_1}|\D(v\theta\sigma)|^2\le
8 |A| \intop_{Q_1}|v|^2
\end{equation}
as ${d\over dt}\intop_{Q_1}v^2\theta^2\sigma^2=\intop_{B_1}v^2(x,0)\theta^2(x)dx\ge 0$, we have from \eqref{2.4.45} and Korn's inequality
\begin{equation}\label{2.4.4}
\intop_{Q_{1/2}}|\nabla v|^2\le C \left[  \intop_{Q_{1/2}}|\D v|^2+\intop_{Q_{1/2}}
|v|^2 \right] \le C \intop_{Q_1}|v|^2
\end{equation}
with constant $C$ depending on $\lambda, |A|, K_p$. By linearity of \eqref{linCQ1} arbitrary derivative of $v$ satisfies again \eqref{linCQ1}. Because the  time derivative can be expressed by space derivative via equation, one can iterate (\ref{2.4.4}) obtaining (with a slight abuse of cutoff function, which should have been stated for $Q_{1-2^{-n}}$)
\begin{equation}\label{2.4.5}
\intop_{Q_{1/2}}|\nabla^{(2k)} v|^2 + |\partial_t^{(k)} v|^2   \le C (\lambda, |A|, K_p, k) |v|_{L^2(Q_1)}^2
\end{equation}
which from Sobolev imbedding, for $k$ big enough, gives
\begin{equation}\label{2.4.6}
\sup_{Q_{1/2}}|\nabla^{(2m)} v|^2 + |\partial_t^{(m)} v|^2   \le C (\lambda, |A|, K_p, m, n) |v|_{L^2(Q_1)}^2
\end{equation}
and as a result for arbitrary $q \ge 1$ and  $\rho/r \le 1/2$ both
\begin{equation}\label{2.4.65}
\left[ \dashint_{Q_{\rho/r}} \left|\nabla^{(2m)} v \right|^q + \left|\partial_t^{(m)} v \right|^q   \right]^\frac{1}{q} \le C (\lambda, |A|, K_p, m, n) |v|_{L^2(Q_1)}
\end{equation}
and
\begin{equation}\label{2.4.9}
\left[ \dashint_{Q_{\rho/r}} \left|v^{(m-1)}-v^{(m-1)}(\tilde{z}_0) \right|^q \right]^\frac{1}{q}\le \left(\frac{\varrho}{r} \right) C (\lambda, |A|, K_p, m, n) |v|_{L^2(Q_1)}
\end{equation}
\begin{equation}\label{2.4.95}
\left[  \dashint_{Q_{\rho/r}} \left|v^{(m-1)}-\left(v^{(m-1)}\right)_{Q_{\rho/r}} \right|^q \right]^\frac{1}{q} \le \left(\frac{\varrho}{r} \right) C (\lambda, |A|, K_p, m, n) |v|_{L^2(Q_1)}^2
\end{equation}
hold, where $\tilde{z}_0$ is an arbitrary point, $v^{(m)}$ denotes $(\nabla^{(2m)} v, \partial_t^{(m)} v)$. To obtain the last two inequalities one uses
\begin{equation}\label{2.4.10}
\left[  \dashint_{Q_{\rho/r}}  |h-h(\tilde{z}_0)|^q \right]^\frac{1}{q} \le  \left(\frac{\varrho}{r} \right)\sup_{Q_{\rho/r}} (|\nabla h|+|h_{,t}|)
\end{equation}
Rescaling inequalities \eqref{2.4.9} --- \eqref{2.4.10} from $v$ back to $u$ we obtain
\begin{equation}\label{2.4.65'}
\left[ \dashint_{Q_{\rho}} \left|\nabla^{(2m)} u \right|^q + \left|\partial_t^{(m)} u \right|^q   \right]^\frac{1}{q} \le r^{-2m} C (\lambda, |A|, K_p, m, n) \left[  \dashint_{Q_{r}} |u|^2 \right]^{1/2}
\end{equation}
\begin{equation}\label{2.4.9'}
\left[ \dashint_{Q_{\rho}} \left|u^{(m-1)}-u^{(m-1)}(\tilde{z}_0) \right|^q \right]^\frac{1}{q}\le  r^{-2(m-1)} \left(\frac{\varrho}{r} \right) C (\lambda, |A|, K_p, m, n) \left[   \dashint_{Q_{r}} |u|^2 \right]^{1/2}
\end{equation}
\begin{equation}\label{2.4.95'}
\left[  \dashint_{Q_{\rho}} \left|u^{(m-1)}-\left(u^{(m-1)}\right)_{Q_{\rho}} \right|^q \right]^\frac{1}{q} \le  r^{-2(m-1)}  \left(\frac{\varrho}{r} \right) C (\lambda, |A|, K_p, m, n) \left[  \dashint_{Q_{r}} |u|^2 \right]^{1/2}
\end{equation}
which together with Giusti's technique, allowing in local inequalities to decrease the power of integrability on right-hand-sides below $2$  by interpolation, implies thesis.
\end{proof}
Let us now turn to the proof of symmetric caloric approximation lemma.
\begin{proof}[Proof of Lemma \ref{lem5.1}]
First consider unit cylinder, i.e. ${Q_{\varrho} (z_0)} = Q_1$. We perform an indirect proof. Contradiction yields existence of  $\varepsilon_0>0$, $\lambda_0>0$, $M_0>0$ such that for any $k\in\N$ one can find such $A_k\in S(\lambda_0, M_0)$,  $\gamma_k \in [0,1] $ that
\begin{equation}\label{cal_c}
f_k\in  H(1; {1/ k}, A_k, \gamma_k) \quad \wedge \quad 
\hforall_{h\in H(1/2;  A_k, \gamma_k)  } \diagintop_{Q_{1/2}}
4 \left|{h-f_k}\right|^2 + \gamma_k^{p-2} 2^p \left|{h-f_k}\right|^p  \ge \varepsilon_0
\end{equation}
By $\gamma_k \in [0,1]$, definitions of approximatively weakly symmetrical caloric functions and $S(\lambda_0, M_0)$ we have the following convergences for $k \rightarrow \infty$ (up to non-relabeled subsequence)
\begin{equation} \label{5.1.15}
\begin{aligned}
\gamma_k \rightarrow& \gamma \quad&\quad
A_k \rightarrow& A \in S(\lambda_0, M_0)  \\
f_k \rightharpoonup& f  \; \text{ in } \;  L^2 (Q_1)   \quad&\quad
\nabla f_k \rightharpoonup& \nabla f \; \text{ in } \; L^2 (Q_1) \\
\gamma_k^{\frac{p-2}{p}} f_k \rightharpoonup& \gamma^{\frac{p-2}{p}}   f  \; \text{ in } \;  L^p (Q_1)   \quad&\quad
 \gamma_k^{\frac{p-2}{p}}  \nabla f_k \rightharpoonup& \gamma^{\frac{p-2}{p}}   \nabla f  \; \text{ in } \;  L^p (Q_1) \\
\end{aligned}
\end{equation}
The fact that $A$ is symmetrizing comes from $A_k$ being symmetrizing. The last two limits can be identified thanks to previous convergences. Consequently by l.w.s.c. of Lebesgue norms one has
\begin{equation} \label{5.1.16}
\diagintop_{Q_1 } \left|{f} \right|^2+ |\nabla f|^2 + \gamma^{p-2} \left[ \left| {f} \right|^p+ |\nabla f|^p \right] \le 1
\end{equation}
Write for arbitrary ${\varphi\in C_0^\infty({Q_1 })}$
\begin{align*}
&\intop_{Q_1} f\varphi_{,t}-A(\D f,\D\varphi)=
\intop_{Q_1} (f-f_k)\varphi_{,t}-\intop_{Q_1} A(\D f-\D f_k,\D\varphi)
\quad-\intop_{Q_1} (A-A_k)(\D f_k,\D\varphi)+\intop_{Q_1} f_k\varphi_{,t}-A_k(\D f_k,\D\varphi)=\\
&I+II+III+IV
\end{align*}
one has that $I,II,III\mathop{\to}\limits^{k\to\infty}0$ by (\ref{5.1.15}), $IV\mathop{\to}\limits^{k\to\infty}0$ as $f_k\in  H(1; {1/ k}, A_k, \gamma_k)$;  therefore
\begin{equation}\label{5.1.2}
\intop_{Q_1} f\varphi_{,t}-A(\D f,\D\varphi)=0 \quad\hforall_{\varphi\in C_0^\infty({Q_1})}
\end{equation}
\eqref{5.1.2} and \eqref{5.1.16} imply 
\begin{equation}\label{5.1.25}
f \in H(1; A, \gamma)
\end{equation}
i.e. limit $f$ is in the set of caloric functions.\\
Assume for a moment that we have also strong convergences
\begin{equation}\label{5.1.4}
f_k\to f\quad {\rm in}\ \ L^2 (Q_1) , \qquad  \gamma_k ^{\frac{p-2}{p}} f_k \to  \gamma^{\frac{p-2}{p}} f \quad {\rm in}\ \ L^p  (Q_1)
\end{equation}
then we would have in view of \eqref{5.1.4}
\begin{equation}\label{cal_ca}
\lim_{k \to \infty} \diagintop_{Q_1}
\left| {f-f_k}\right|^2 + \gamma_k^{p-2} \left| {f-f_k}\right|^p  = 0
\end{equation}
which with \eqref{5.1.25} is almost a contradiction to \eqref{cal_c} with an exception, that we require in \eqref{cal_c}  $h \in H(1/2;  A_k, \gamma_k) $ and instead have $f \in H(1; A, \gamma)$. We compensate this difference by proceeding as follows. Consider the following linear boundary-value problem 
\begin{align}\label{pom}
\begin{cases}
\omega^k_{,t}-\divv A_k\D\omega^k=0 \quad &{\rm in}\ \ Q_{3/4}\\
\omega^k=f \quad &\text{on } \partial_\Gamma Q_{3/4}
\end{cases}
\end{align}
 \eqref{pom} and \eqref{5.1.2} give
\begin{equation}\label{5.1.5}
\intop_{Q_{3/4}} (\omega^k -f) \varphi_{,t}-\intop_{Q_{3/4}} A_k(\D\omega^k - \D f,\D\varphi)=\intop_{Q_{3/4}} (A_k-A)(\D f,\D\varphi)
\end{equation}
since by \eqref{pom} $w^k$ and $f$ agree on  the boundary, we can test \eqref{5.1.5} with $w^k - f$ obtaining
\begin{equation*}
\sup_{t \in (- (3/4)^2, 0)} |\omega^k(t) - f(t)|_{L^2(B_{3/4})}^2+\intop_{B_{3/4}} A_k(\D (\omega^k -f),\D(\omega^k -f))=\intop_{B_{3/4}}  (A_k-A)(\D f,\D (\omega^k -f))
\end{equation*}
which by ellipticity of $A_k$ and Korn's inequality gives
\begin{equation}\label{5.1.6}
|\omega^k -f|_{L^2({Q_{3/4}}) }+|\nabla (\omega^k - f)|_{L^2({Q_{3/4}}) } \le C \left[\sup_t|\omega^k(t) - f(t)|_{L^2({B_{3/4}})}+|\D (\omega^k-f)|_{L^2({Q_{3/4}})} \right]\le
C|A_k-A|\mathop{\to}\limits^{k\to\infty}0
\end{equation}
Observe that in view of inequality \eqref{2.4.1} of Lemma \ref{LinSys}, $|\omega^k|_{L^2 (Q_{3/4})} $ controls norms on $ Q_{1/2} $ of $\omega^k$ of arbitrary high order. This and \eqref{5.1.6} yield
\begin{equation}\label{5.1.65}
| (\omega^k - f)|_{L^p (Q_{1/2})} + |\nabla (\omega^k - f)|_{L^p (Q_{1/2})} \to 0
\end{equation}
We show now that  $\omega_k$ contradicts \eqref{cal_c} for large $k$. One has 
\begin{equation}\label{5.1.66}
\dashint_{Q_{1/2}} \gamma_k^{p-2} |\omega_k-f_k|^p \le \dashint_{Q_{1/2}}|  \gamma_k^{p-2} f- \gamma_k^{p-2}  f_k|^p  + \dashint_{Q_{1/2}} |\omega_k-f|^p  \to 0
\end{equation}
where the convergence stems from \eqref{5.1.65} and \eqref{5.1.4}. Similarly
\begin{equation}\label{5.1.67}
\dashint_{Q_{1/2}}|\omega_k-f_k|^2 \to 0
\end{equation}
What's more, $\gamma_k \rightarrow \gamma$, \eqref{5.1.6} and \eqref{5.1.65} imply
\begin{equation}\label{5.1.7}
\dashint_{Q_{1/2}}|\omega_k|^2 + |\nabla \omega_k|^2 + \gamma_k^{p-2} [|\omega_k|^p + |\nabla \omega_k|^p] \to 
\dashint_{Q_{1/2}}|f|^2 + |\nabla f|^2 + \gamma^{p-2} [|f|^p + |\nabla f|^p] \le 2^{n+2}
\end{equation}
where the last inequality holds thanks to $f \in H(1; A, \gamma)$. \eqref{5.1.7} and \eqref{pom} state that $\omega_k \in H(1/2;  A_k, \gamma_k)  $ for large $k$, whence \eqref{5.1.66}, \eqref{5.1.67} show that $\omega_k$ approximates $f_k$ in $L^2 - L^p$ sense. This contradicts \eqref{cal_c}. Thus we are done with the case of an unit parabolic cylinder, provided (\ref{5.1.4}) really holds, which now will be proven by means of parabolic compactness of  Lemma \ref{lem1.15}. By \eqref{cal_c} $  f_k\in  H(1; {1/ k}, A_k, \gamma_k) $, so from its definition
\begin{equation*}
\left|\intop_Q \gamma_k^{\frac{p-2}{p}} f_k\varphi_{,t}\right|\le {\lambda^{-1}_0}\left(\intop_Q|\D\varphi|^{p'}\right)^{1\over {p'}}+{1\over k}\sup_Q|\D\varphi|
\end{equation*}
taking $\varphi=\xi(x)\eta_{s_1,s_2}^\varepsilon(t)$ with
\begin{equation*}
\eta=\left\{\begin{aligned}
&1\quad &t\in[s_1,s_2]\\
&0\quad &t\not\in[s_1-\varepsilon,s_2+\varepsilon]\\
&{\rm affine\ otherwise}\end{aligned}\right.
\end{equation*}
\begin{equation*}
\left|\intop_B{1\over\varepsilon}\left[\intop_{s_1-\varepsilon}^{s_1} \gamma_k^{\frac{p-2}{p}} f_k-
\intop_{s_2}^{s_2+\varepsilon} \gamma_k^{\frac{p-2}{p}} f_k\right]\xi\right|\le\left[ \lambda^{-1}_0 {(s_2-s_1+2\varepsilon)}^{\frac{1}{p'}}+
{1\over k}\right]|\D\xi|_{L^\infty(B)}
\end{equation*}
so for a.e. $s_1,s_2$, by $\varepsilon\to0$,
\begin{equation}\label{5.1.9}
\left|\intop_B( \gamma_k^{\frac{p-2}{p}}  f_k(\cdot,s_2)-  \gamma_k^{\frac{p-2}{p}} f_k(\cdot,s_1))\xi\right|\le\left(c (s_2-s_1)^\frac{1}{p'} +
{1\over k}\right)|\nabla\xi|_{L^\infty(B)}
\end{equation}
for $l-{n\over2}>1$ holds $C|\xi|_{W^{l,2}}\ge|\nabla\xi|_{L^\infty}$ so (\ref{5.1.9}) and a density argument imply
\begin{equation*}
|  \gamma_k^{\frac{p-2}{p}}  f_k(\cdot,s_2)-  \gamma_k^{\frac{p-2}{p}}  f_k(\cdot,s_1)|_{W^{-l,2}(B)}\le c\left(  (s_2-s_1)^\frac{1}{p'} +{1\over k}\right)
\end{equation*}
i.e.
\begin{equation}\label{5.1.10}
\intop_{-1}^{-h}| \gamma_k^{\frac{p-2}{p}}  f_k(\cdot,t+h)-  \gamma_k^{\frac{p-2}{p}}  f_k(\cdot,t)|^p_{W^{-l,2}(B)}\le c\left(h^{p-1}+{1\over k^{p'}}\right)
\end{equation}
Fix $\varepsilon > 0$. For $k > (2c/ \varepsilon)^\frac{1}{p'} $ r.h.s. of \eqref{5.1.10} does not exceed $\varepsilon $ for every $h \le  ( \varepsilon / 2c)^\frac{1}{p - 1}$. For a finite number of initial $k$ one can majorize l.h.s. of \eqref{5.1.10} by $\varepsilon $  for every $h \le h_0$ by properties of Bochner spaces. All in all it holds
\begin{equation}\label{5.1.11}
\hforall_{\varepsilon > 0} \hexists_{h'} \hforall_{h \le h'} \intop_{-1}^{-h}| \gamma_k^{\frac{p-2}{p}}  f_k(\cdot,t+h)-  \gamma_k^{\frac{p-2}{p}}  f_k(\cdot,t)|^p_{W^{-l,2}(B)}\le \varepsilon
\end{equation}
This and uniform boundedness of $\gamma_k^{\frac{p-2}{p}}  f_k$ in $L^p (W^{1,p})$ allows us to use Lemma \ref{lem1.15} with choices $X = W^{1,p}, \; Y = L^p, \; Z = W^{-l,2}$ to get, up to subsequence, $\gamma_k^{\frac{p-2}{p}}  f_k \to \gamma^{\frac{p-2}{p}}  f$ in $L^p$. Similar argument for $L^2$ gives  (\ref{5.1.4}). \\
Finally, we perform step from the unit cylinder $Q_1 (z_0)$ to an arbitrary one $Q_\varrho(z_0)$. To this end it suffices to observe that mapping $f (\cdot) \rightarrow \frac{1}{r} f ( \frac{\cdot}{r})$ is a bijection from $H (r; \delta, A, \gamma)$ to $H (1; \delta, A, \gamma)$ and from $H (r; A, \gamma)$ to $H (1; A, \gamma)$.
\end{proof}
It is worth remarking that it seems possible to provide a more constructive proof, based on ideas of Diening and collaborators (see \cite{[DSV]} and references therein).

Below is presented proof of Lemma \ref{lem3.05}, which states some algebraic inequalities useful for local estimates of Section \ref{sec:ineq}.
\begin{proof}[Proof of Lemma \ref{lem3.05}]
Observe first that assumption \eqref{1.8} yields
\begin{equation}\label{1.8'}
|A(z,u,q)-A(\tilde z,\tilde u,  q)|\le C K(2|u|+1)( d_p ( z - \tilde z) +|u-\tilde u|)^\beta(1+|q|^{p-1})
\end{equation}
this with  $|l(z_0)| + |\nabla l| \equiv \tilde M$  gives
\begin{equation}
\begin{aligned}
&|A(z,u,P)-A(z_0,l(z_0),P)| \le C K(2|l(z_0)|+1)  (|z-z_0|+|u-l(z_0)|)^\beta (1+|P|^{p-1})\\
&\le  C(\tilde M) (\varrho^\beta+|u-l|^\beta+|l-l(z_0)|^\beta)
(1+|P-\D l|^{p-1}+|\D l|^{p-1})\\
&\le C(\tilde M) (\varrho^\beta+|u-l|^\beta+ (|\nabla l||z-z_0|)^\beta)(1+|P-\D l |^{p-1} + \tilde M^{p-1})\\
&\le C(\tilde M) \varrho^\beta
\left [1+  |P-\D l|^{p-1} + \left| \frac{u-l}{ \varrho} \right|^{p \beta} + |P - \D l|^{p} + \left| \frac{u - l}{\varrho} \right|^\beta \right]
\le C(\tilde M) \varrho^\beta
\left [1+  |P-\D l|^{p} + \left| \frac{u - l}{\varrho} \right|^p \right] 
\end{aligned}
\end{equation}
i.e. \eqref{alg_1}.\\
Turning to estimate for \eqref{alg_3} one has again from \eqref{1.8'}
\begin{multline}
|A(z_0,l(z_0), \D l) - A(z, u, \D l)| 
\le  C K(| l(z_0)|)(\varrho+|u-l_0|)^\beta(1+|\D l|)^{p-1} \\
\le C(\tilde M) [\varrho^\beta + |u-l|^\beta+
|l-l(z_0)|^\beta] \le C(\tilde M) \varrho^\beta \left[1+\left| \frac{u - l}{\varrho} \right|^\beta \right] 
\end{multline}
To show inequality \eqref{alg_2} write
\begin{equation}
|A(z,u,P)-A(z,u,\D l)| \le |A(z,u,P)-A(z,l,P)| + |A(z,l,P)-A(z,l,\D l)| +  |A(z,u,\D l)-A(z,l,\D l)| = I + II + III.
\end{equation}
Estimate I using \eqref{1.8} and properties of function $K$
\begin{equation}
\begin{aligned}
& |A(z,u,P)-A(z,l,P)| \le C \cdot \min \left(1, K(|u|+|l|)   |u-l|^\beta \right)  \cdot (1+|P|^{p-1}) \le \\
& C  \min \left(1, K(|u-l|+2|l|)  \cdot |u-l|^\beta \right)  \cdot (1+ |P-\D l|^{p-1}+|\D l|^{p-1}) \le  \\
&C (\tilde M)  \min \left(1, K(|u-l|+2|l-l(z_0)|+ 2\tilde M)  \cdot |u-l|^\beta \right)  \cdot (1+ |P-\D l|^{p-1}) \end{aligned}
\end{equation}
When $|u-l| \ge 1$ we estimate further by $C(\tilde M) (1+ |P-\D l|^{p-1}) \le C(\tilde M) (|u-l| + |P-\D l|^{p-1})$.\\
Otherwise by 
$C (\tilde M)  \min \left(1, K(1+2\varrho |\nabla l|+ 2\tilde M)  \cdot |u-l|^\beta \right)  \cdot (1+ |P-\D l|^{p-1})  \le C (\tilde M)  \left( |u-l|^\beta+ |P-\D l|^{p-1}  \right)  $.\\
Consequently
\begin{equation}
I \le C (\tilde M)  \left( |u-l|^\beta+ |u-l| + |P-\D l|^{p-1}  \right) 
\end{equation}
Next consider $II$ in case $| P - \D l | \le 1$. In view of $|l| \le \tilde M+ |\nabla l| \varrho \le 2\tilde M$ we estimate $II$ by \eqref{1.7} as follows
\begin{equation}
II = \left| \int_0^1 \frac{\partial A}{\partial q}  (z, l, s ( P - \D l ) + \D l) ( P - \D l ) ds \right| \le  C_\eqref{1.7} (3\tilde M+1) | P - \D l |
\end{equation}
where in case $| P - \D l | \ge 1$ we robustly estimate  $II$  by \eqref{1.6}
\begin{equation}
II  \le C \cdot (1+|\D l|^{p-1} + |P|^{p-1}) \le C(\tilde M) (1+|\D l - P|^{p-1}) \le  C(\tilde M) (|\D l - P|+|\D l - P|^{p-1})
\end{equation}
so from both cases one has
\begin{equation}
II \le C(\tilde M) (|\D l - P|+|\D l - P|^{p-1})
\end{equation}
Finally we obtain estimate for $III$ in view of \eqref{1.8'}
\begin{equation}
III \le C(\tilde M) |u-l|^\beta
\end{equation}
and inequalities for $I, II, III$ give desired \eqref{alg_2}. Inequality  \eqref{alg_4} follows similarly and more straightforwardly, when instead of \eqref{1.8'} we observe that assumption  \eqref{1.8} implies
\begin{equation}
|A(z,u,q)-A(\tilde z,\tilde u,  q)|\le  C \cdot (1+|q|^{p-1}).
\end{equation}
\end{proof}
\vskip 1cm


\end{document}